\documentclass[12pt,a4paper,twoside]{article}

\usepackage[latin1]{inputenc}
\usepackage[T1]{fontenc}
\usepackage{graphicx}
\usepackage{latexsym}
\usepackage{amsmath,amssymb,amsthm}
\usepackage{bbm}
\usepackage{ifpdf}
\usepackage{exscale}
\usepackage{fancyhdr}
\usepackage{tocloft}

\usepackage[DIV12]{typearea}
\newcommand{\e}{\ensuremath{\varepsilon}}
\newcommand{\dt}{\ensuremath{\textup{d}}}
\newcommand{\dist}{\operatorname{d}}
\newcommand{\dL}{\operatorname{d}_L}

\newcommand{\dtL}{\operatorname{\tilde{d}}}
\newcommand{\sh}{\operatorname{Sh}}

\newcommand{\czero}{\operatorname{\bf{A0}}}
\newcommand{\cone}{\operatorname{\bf{A1}}}
\newcommand{\cB}{\operatorname{\bf{B}}}
\newcommand{\cspace}{\operatorname{\bf{C1}}}
\newcommand{\ctime}{\operatorname{\bf{C2}}}

\newcommand{\good}{\operatorname{Good}_L}
\newcommand{\badP}{\ensuremath{\mathcal{B}}}
\newcommand{\lambdap}{\ensuremath{\lambda^{(p)}}}
\newcommand{\goods}{\operatorname{Good}_L^{\textup{\tiny sp}}}
\newcommand{\bads}{\operatorname{Bad}_L^{\textup{\tiny sp}}}
\newcommand{\goodt}{\operatorname{Good}_L^{\textup{\tiny tm}}}
\newcommand{\onebadt}{\operatorname{OneBad}_L^{\textup{\tiny tm}}}
\newcommand{\manybadt}{\operatorname{ManyBad}_L^{\textup{\tiny tm}}}
\newcommand{\badPs}{\ensuremath{\mathcal{B}}_L^{\textup{\tiny sp}}}
\newcommand{\badPt}{\ensuremath{\mathcal{B}}_L^{\textup{\tiny tm}}}

\newcommand{\pip}{\ensuremath{\pi^{(p)}}}
\newcommand{\piq}{\ensuremath{\pi^{(q)}}}
\newcommand{\ph}{\ensuremath{\hat{\pi}}}
\newcommand{\php}{\ensuremath{\hat{\pi}^{(p)}}}
\newcommand{\phq}{\ensuremath{\hat{\pi}^{(q)}}}
\newcommand{\ptp}{\ensuremath{\tilde{\pi}^{(p)}}}

\newcommand{\Ph}{\ensuremath{\hat{\Pi}}}
\newcommand{\Phg}{\ensuremath{\hat{\Pi}^g}}

\newcommand{\gh}{\ensuremath{\hat{g}}}
\newcommand{\ghp}{\ensuremath{\hat{g}}^{(p)}}
\newcommand{\ghq}{\ensuremath{\hat{g}}^{(q)}}
\newcommand{\gtp}{\ensuremath{\tilde{g}}^{(p)}}

\newcommand{\Gh}{\ensuremath{\hat{G}}}
\newcommand{\Ghg}{\ensuremath{\hat{G}^g}}

\newcommand{\Prw}{\operatorname{P}}
\newcommand{\Erw}{\operatorname{E}}
\newcommand{\Prwt}{\ensuremath{\widetilde{\Prw}}}
\newcommand{\Erwt}{\ensuremath{\widetilde{\Erw}}}
\newcommand{\pP}{\ensuremath{\mathbb{P}}}
\newcommand{\pQ}{\ensuremath{\mathbb{Q}}}
\newcommand{\pE}{\ensuremath{\mathbb{E}}}
\newcommand{\nP}{\ensuremath{\textup{P}}}
\newcommand{\nE}{\ensuremath{\textup{E}}}
\newcommand{\pF}{\ensuremath{\mathcal{F}}}
\newcommand{\pG}{\ensuremath{\mathcal{G}}}
\newcommand{\mP}{\ensuremath{\mathcal{P}}}
\newcommand{\mPsymone}{\ensuremath{\mathcal{P}^{\textup{s},1}}}
\newcommand{\mPsym}{\ensuremath{\mathcal{P}^{\textup{s}}_\kappa}}
\newcommand{\mPsymi}{\ensuremath{\mathcal{P}^{\textup{s}}_\iota}}
\newcommand{\bj}{\ensuremath{{\bf j}}}
\newcommand{\mJ}{\ensuremath{\mathcal{J}_m}}

\newcommand{\df}{\ensuremath{\overset{\mbox{\scriptsize\textup{def}}}{=}}}

\newtheorem{theorem}{Theorem}[section]
\newtheorem{lemma}{Lemma}[section]
\newtheorem{proposition}{Proposition}[section]
\newtheorem{corollary}{Corollary}[section]
\newenvironment{proof1}{\noindent{\bf Proof:}}{%
  \hspace*{\fill}$\Box$\par\vskip2ex}
\newenvironment{proof2}{\noindent{\bf Proof }}{%
  \hspace*{\fill}$\Box$\par\vskip2ex}

\theoremstyle{definition}
\newtheorem{remark}{Remark}[section]

\begin{document}
\pagestyle{empty} 
\pagestyle{fancy}
\renewcommand{\headrulewidth}{0mm} \chead[\leftmark]{\rightmark}
\rhead[]{\thepage} \lhead[\thepage]{} \cfoot{} 
\pagenumbering{arabic}

\title{An invariance principle for a class of non-ballistic random walks in random
  environment} 

\author{Erich Baur\footnote{Email:
    erich.baur@math.uzh.ch.}\\ENS Lyon}

\maketitle 
\thispagestyle{empty}
\begin{abstract}
  We are concerned with random walks on $\mathbb{Z}^d$, $d\geq 3$, in an
  i.i.d. random environment with transition probabilities $\e$-close to
  those of simple random walk. We assume that the environment is balanced
  in one fixed coordinate direction, and invariant under reflection in the
  coordinate hyperplanes. The invariance condition was used in~\cite{BB} as
  a weaker replacement of isotropy to study exit distributions. We obtain
  precise results on mean sojourn times in large balls and prove a quenched
  invariance principle, showing that for almost all environments, the
  random walk converges under diffusive rescaling to a Brownian motion with
  a deterministic (diagonal) diffusion matrix. We also give a concrete
  description of the diffusion matrix. Our work extends the results of
  Lawler~\cite{LAWbalanced}, where it is assumed that the environment is
  balanced in all coordinate directions.
\end{abstract}

   \let\thefootnote\relax
  \footnotetext{MSC 2000 {\it subject classifications.}
    Primary 60K37; secondary 82C41.}
  \footnotetext{{\it Acknowledgment of support.} This research was
    supported by the Swiss National Science Foundation grant
    P2ZHP2\_151640.}

\section{Introduction and main results}
\subsection{The model}
Denote by $e_i$ the $i$th unit vector of $\mathbb{Z}^d$. We let
$\mathcal{P}$ be the set of probability distributions on $\{\pm e_i : i =
1,\dots,d\}$  and put $\Omega=\mathcal{P}^{\mathbb{Z}^d}$. Denote by $\pF$
the natural product $\sigma$-field
on $\Omega$ and by $\pP = \mu^{\otimes\mathbb{Z}^d}$  the product
probability measure on $(\Omega,\pF)$. 

Given an element (or environment) $\omega\in\Omega$, we denote by
$(X_n)_{n\geq 0}$ the canonical nearest neighbor Markov
chain on $\mathbb{Z}^d$ with transition probabilities
$$
  p_\omega(x,x+e) = \omega_x(e),\quad e\in\{\pm e_i : i = 1,\dots,d\},
$$
the {\it random walk in random environment} (RWRE
for short). We write $\Prw_{x,\omega}$ for the ``quenched'' law of $(X_n)_{n\geq 0}$
starting at $x\in\mathbb{Z}^d$.

We are concerned with RWRE in dimensions $d\geq 3$ which is an
$\e$-perturbation of simple random walk. To fix a perturbative regime, we
shall assume the following condition.
\begin{itemize}
\item Let $0<\varepsilon<1/(2d)$. We say that $\czero(\e)$ holds if
  $\mu(\mathcal{P}_\e) = 1$, where
$$
\mathcal{P}_\e = \left\{q\in\mathcal{P} : \left|q(\pm e_i) -
  1/(2d)\right| \leq \e\mbox{ for all }i=1,\dots,d\right\}. 
$$
\end{itemize} 
Furthermore, we work under two centering conditions on the measure
$\mu$. The first rules out ballistic behavior, while the second
guarantees that the RWRE is balanced in direction $e_1$.
\begin{itemize}
\item We say that $\cone$ holds if 
  $\mu$ is invariant under all $d$ reflections
  $O_i:\mathbb{R}^d\rightarrow\mathbb{R}^d$ mapping the unit vector $e_i$ to
  its inverse, i.e. $O_ie_i=-e_i$ and $O_ie_j=e_j$ for $j\neq i$.
  In other words, the laws of $(\omega_0(O_ie))_{|e|=1}$ and
  $(\omega_0(e))_{|e|=1}$ coincide, for each $i=1,\dots,d$.
\item We say that $\cB$ holds if $\mu(\mPsymone)=1$, where
$$\mPsymone = \{p\in\mathcal{P} : p(e_1)=p(-e_1)\}.$$
\end{itemize}
We now state our results. Then we discuss them together with our conditions
in the context of known results from the literature.

\subsection{Our main results}
Our first statement shows $\pP$-almost sure convergence of the (normalized)
RWRE mean sojourn time in a ball when its radius gets larger and
larger. Let $V_L=\{y\in\mathbb{Z}^d:|y|\leq L\}$ denote the discrete ball
of radius $L$, and $V_L(x) = x +V_L$. Denote by $\tau_{V_L(x)}=\inf\{n\geq
0:X_n\notin V_L(x)\}$ the first exit time of the RWRE from $V_L(x)$. We
write $\Erw_{x,\omega}$ for the expectation with respect to
$\Prw_{x,\omega}$.  We always assume $d\geq 3$.

\begin{theorem}[Quenched mean sojourn times, $d\geq 3$]
\label{thm-times}
Assume $\cone$ and $\cB$. Given $0<\eta < 1$, one can find
$\e_0=\e_0(\eta) > 0$ such that if $\czero(\e)$ is satisfied for some
$\e\leq \e_0$, then the following holds: There exists $D\in [1-\eta,
1+\eta]$ such that for $\pP$-almost all $\omega\in\Omega$,
$$
\lim_{L\rightarrow\infty}\left(\Erw_{0,\omega}\left[\tau_{V_L}\right]/L^2\right)
= D.
$$
Moreover, one has for each $k\in\mathbb{N}$, for $\pP$-almost all $\omega$,
$$\lim_{L\rightarrow\infty}\left(\inf_{x: |x| \leq
    L^k}\Erw_{x,\omega}\left[\tau_{V_L(x)}\right]/L^2 \right)=
\lim_{L\rightarrow\infty}\left(\sup_{x: |x| \leq
    L^k}\Erw_{x,\omega}\left[\tau_{V_L(x)}\right]/L^2 \right)= D. 
$$
\end{theorem}

Standard arguments then imply the following bound on the moments.
\begin{corollary}[Quenched moments]
 \label{cor-times-moments}
 In the setting of Theorem~\ref{thm-times}, for each $k$,
 $m\in\mathbb{N}$ and $\pP$-almost all $\omega$,
$$
\limsup_{L\rightarrow\infty}\left(\sup_{x: |x| \leq
    L^k}\Erw_{x,\omega}\left[\tau^m_{V_L(x)}\right]/L^{2m} \right)\leq
2^mm!\,.
$$
\end{corollary}

Combined with results on the spatial behavior of the RWRE from~\cite{BB},
we prove a functional central limit theorem under the quenched measure. 
In~\cite{BB}, it was shown that under $\cone$ and $\czero(\e)$ for $\e$
small, the limit 
\begin{equation}
\label{eq-pinfty}
2p_\infty(\pm e_i) =
\lim_{L\rightarrow\infty}\sum_{y\in\mathbb{Z}^d}\pE\left[\Prw_{0,\omega}\left(X_{\tau_{V_L}}=y\right)\right]\frac{y_i^2}{|y|^2} 
\end{equation}
exists for $i=1,\dots,d$, and $|p_\infty(e_i)-1/(2d)|\rightarrow 0$ as
$\e\downarrow 0$. Let
\begin{equation}
\label{covariance-matrix-def}
{\bf \Lambda}=\left(2\,p_\infty(e_i)\delta_{i}(j)\right)_{i,j=1}^d\in\mathbb{R}^{d\times d},
\end{equation}
and define the linear interpolation  
$$
X_t^{n}= X_{\lfloor tn\rfloor} + (tn-\lfloor
  tn\rfloor)\left(X_{\lfloor tn\rfloor+1}-X_{\lfloor
    tn\rfloor}\right),\quad t\geq 0.
$$
The sequence $(X_t^{n}, t\geq 0)$ takes values in the space
$C(\mathbb{R}_+,\mathbb{R}^d)$ of
$\mathbb{R}^d$-valued continuous functions on $\mathbb{R}_+$. The set $C(\mathbb{R}_+,\mathbb{R}^d)$ is
tacitly endowed with the uniform topology and its Borel $\sigma$-field.

\begin{theorem}[Quenched invariance principle, $d\geq 3$]
\label{thm-clt}
Assume $\czero(\e)$ for small $\e > 0$, $\cone$ and $\cB$. Then for
$\pP$-a.e. $\omega\in\Omega$, under $\Prw_{0,\omega}$, 
\begin{quote}
$X_\cdot^{n}/\sqrt{n}$ converges in
law to a $d$-dimensional Brownian motion with diffusion matrix
$D^{-1}{\bf \Lambda}$, where $D$ is the constant from Theorem~\ref{thm-times} and
${\bf \Lambda}$ is given by~\eqref{covariance-matrix-def}.
\end{quote}
\end{theorem}

Under Conditions $\czero(\e)$ and $\cone$, a local limit law for RWRE exit
measures was proved in~\cite{BB} for dimensions three and higher. Before
that, similar results were obtained by Bolthausen and Zeitouni~\cite{BZ}
for the case of isotropic perturbative RWRE. While the results of~\cite{BZ}
and~\cite{BB} do already imply transience for the random walks under
consideration, they do not prove diffusive behavior, since there is no
control over time. This was already mentioned in~\cite{BZ}: ``In future
work we hope to combine our exit law approach with suitable exit time
estimates in order to deduce a (quenched) CLT for the RWRE.''  Under the
additional Condition $\cB$ which shall be discussed below, we fulfill here
their hope.

In dimensions $d>1$, the RWRE under the quenched measure is an irreversible
(inhomogeneous) Markov chain. A major difficulty in its analysis comes from
the presence of so-called traps, i.e. regions where the random walk can
hardly escape and therefore spends a lot of time. In the ballistic regime
where the limit velocity $\lim_{n\rightarrow\infty}X_n/n$ is an almost sure
constant different from zero, powerful methods leading to law of large
numbers and limit theorems have been established, see
e.g. Sznitman~\cite{SZ1,SZ2,SZ4}, Berger~\cite{BER} or the lecture notes of
Sznitman~\cite{SZ-LN2} with further references. They involve the
construction of certain regeneration times, where, roughly speaking, the
walker does not move ``backward'' anymore.

In the non-ballistic case, different techniques based on renormalization
schemes are required.  In the small disorder regime, results can be found
under the classical isotropy condition on $\mu$, which is stronger than our
condition $\cone$. It requires that for any orthogonal map $\mathcal{O}$
acting on $\mathbb{R}^d$ which fixes the lattice $\mathbb{Z}^d$, the law of
$(\omega_0(\mathcal{O}e))_{|e|=1}$ and $(\omega_0(e))_{|e|=1}$
coincide. Under this condition, Bricmont and Kupiainen~\cite{BK} provide a
functional central limit theorem under the quenched measure for dimensions
$d\geq 3$. However, it is of a certain interest to find a new
self-contained proof of their result. A continuous counterpart is studied
in Sznitman and Zeitouni~\cite{SZZ}. For $d\geq 3$, they prove a quenched
invariance principle for diffusions in a random environment which are small
isotropic perturbations of Brownian motion. Invariance under all lattice
isometries is also assumed in the aforementioned work of Bolthausen and
Zeitouni~\cite{BZ}.
 
In the non-perturbative setting, Bolthausen {\it et al.}~\cite{BSZ} use
so-called cut times as a replacement of regeneration times. At such 
times, past and future of the path do not intersect. However, in
order to ensure that there are infinitely many cut times, it is
assumed in~\cite{BSZ} that the projection of the RWRE onto at least
$d_1\geq 5$ components behaves as a standard random walk. Among other
things, a quenched invariance principle is proved when $d_1\geq 7$ and the
law of the environment is invariant under the antipodal transformation
sending the unit vectors to their inverses.

Our Condition $\cB$ requires only that the environment is balanced in one
fixed coordinate direction ($e_1$ for definiteness). Then the projection of
the RWRE onto the $e_1$-axis is a martingale under the quenched measure,
which implies a priori bounds on the sojourn times, see the
discussion in Section~\ref{Stimes}.  Clearly, assuming just Conditions
$\czero(\e)$ and $\cB$ could still result in ballistic behavior, but the
combination of $\czero(\e)$, $\cone$ and $\cB$ provides a natural framework
to investigate non-ballistic behavior of ``partly-balanced'' RWRE in the
perturbative regime.

To our knowledge, we are the first who study random walks in random
environment which is balanced in only one coordinate direction. The study
of fully balanced RWRE when $\pP(\omega_0(e_i)=\omega_0(-e_i)\ \hbox{for
  all } i=1,\dots,d)=1$ goes back to Lawler~\cite{LAWbalanced}. He proves a
quenched invariance principle for ergodic and elliptic environments in all
dimensions. Extensions within the
i.i.d. setting to the mere elliptic case were obtained by Guo and
Zeitouni~\cite{GZ}, and recently to the non-elliptic case by Berger and
Deuschel~\cite{BD}.

Since the results of~\cite{BB} do also provide local estimates, we believe
that with some more effort, Theorem~\ref{thm-clt} could be improved to a
local central limit theorem. Furthermore, we expect that our results remain
true without assuming Condition $\cB$. Getting rid of this condition would
however require a complete control over large sojourn times, which remains
a major open problem.

\subsubsection{Organization of the paper and rough strategy of the proofs}
We first introduce the most important notation. For ease of readability, we
recapitulate in Section~\ref{SBB} those concepts and results from~\cite{BB}
which play a major role here. In Section~\ref{Scontrolgreen} we provide the
necessary control over Green's functions. To a large extend, we can rely on
the results from~\cite{BB}, but we need additional difference estimates for
our results on mean sojourn times.

In Section~\ref{Stimes}, we prove Theorem~\ref{thm-times}. In this regard,
we shall first show that with high probability, the quenched mean times
$\Erw_{0,\omega}[\tau_L]/L^2$ lie for large $L$ in a small interval
$[1-\eta,1+\eta]$ around $1$ 
(Proposition~\ref{main-prop-times}). This involves the propagation of a
technical Condition $\ctime$ (see Section~\ref{SUBSctimes}).  Once we have
established this, we prove convergence of the (non-random) sequence
$\pE[\Erw_{0,\omega}[\tau_L]]/L^2$ towards a constant $D\in[1-\eta,1+\eta]$
(Proposition~\ref{prop-times-underP}), where $\pE$ denotes the expectation
with respect to $\pP$. Finally, a concentration argument shows that with
high probability, $\Erw_{0,\omega}[\tau_L]$ is close to its mean
$\pE[\Erw_{0,\omega}[\tau_L]]$ (Lemma~\ref{Ltimes-concentration}). This will
allow us to deduce Theorem~\ref{thm-times}.

In the last part of this paper starting with Section~\ref{SCLT}, we show
how Theorem~\ref{thm-times} can be combined with the results on exit laws
from~\cite{BB} to obtain Theorem~\ref{thm-clt}. A strategy of proof of this
statement can be found at the beginning of Section~\ref{SCLT}.

\subsection{Some notation}
We collect here some notation that is frequently used in this
text. 
\subsubsection{Sets and distances}
We put $\mathbb{N} = \mathbb{N}_0= \{0,1,2,3,\dots\}$. For $x\in\mathbb{R}^d$, $|x|$ is the
Euclidean norm of $x$. The distance between $A, B \subset \mathbb{R}^d$ is
denoted $\dist(A,B) = \inf\{|x-y| : x\in A,\; y\in B\}$.  Given $L > 0$, we
let $V_L=\{x\in\mathbb{Z}^d : |x| \leq L\}$, and for $x\in\mathbb{Z}^d$,
$V_L(x) = x + V_L$. Similarly, put $C_L=\{x\in \mathbb{R}^d : |x| <
L\}$. The outer boundary of $V\subset \mathbb{Z}^d$ is given by $\partial V
=\{x\in \mathbb{Z}^d\backslash V: \dist(\{x\},V) = 1\}$. For
$x\in\overline{C}_L$, we let $\dL(x) = L-|x|$. Note for $x\in V_L$,
$\dL(x)\leq \dist(\{x\},\partial V_L)$. Finally, for $0\leq a<b\leq L$, put
$$
\sh_L(a,b)=\{x\in V_L : a\leq \dL(x) < b\},\quad\sh_L(b) = \sh_L(0,b).
$$

\subsubsection{Functions}
We use the usual notation $a\wedge b = \min\{a,b\}$ for reals $a,b$. We
further write $\log$ for the logarithm to the base {\rm e}. 
Given two functions $F,G :
\mathbb{Z}^d\times\mathbb{Z}^d\rightarrow\mathbb{R}$, we write $FG$ for the
(matrix) product $FG(x,y) = \sum_{u\in\mathbb{Z}^d}F(x,u)G(u,y)$, provided
the right hand side is absolutely summable. $F^k$ is the $k$th power
defined in this way, and $F^0(x,y) = \delta_x(y)$. $F$ can also operate on
functions $f:\mathbb{Z}^d\rightarrow\mathbb{R}$ from the left via $Ff(x) =
\sum_{y\in\mathbb{Z}^d}F(x,y)f(y)$.

As usual, $1_W$ stands for the indicator function of the set $W$, but we
will also write $1_W$ for the kernel $(x,y)\mapsto 1_W(x)\delta_x(y)$,
where the Delta function $\delta_x(y)$ is equal to one if $y=x$ and zero
otherwise. If $f:\mathbb{Z}^d\rightarrow \mathbb{R}$, $\|f\|_1 =
\sum_{x\in\mathbb{Z}^d}|f(x)| \in [0,\infty]$ denotes its $L^1$-norm. For a
(signed) measure $\nu : \mathbb{Z}^d\rightarrow \mathbb{R}$, we write
$\|\nu\|_1$ for its total variation norm.

\subsubsection{Transition kernels, exit times and exit measures}
Denote by $\pG$ the $\sigma$-algebra on $(\mathbb{Z}^d)^{\mathbb{N}}$
generated by the cylinder functions. If $p = {(p(x,y))}_{x,y\in \mathbb{Z}^
d}$ is a
family of (not necessarily nearest neighbor) transition probabilities, we
write $\Prw_{x,p}$ for the law of the canonical random walk ${(X_n)}_{n\geq
  0}$ on $({(\mathbb{Z}^d)}^{\mathbb{N}},\pG)$ started from $X_0=x$
$\Prw_{x,p}$ -a.s. and evolving according to the kernel $p$. 

The simple random walk kernel is denoted $p_o(x,x\pm e_i) = 1/(2d)$, and we
write $\Prw_{x}$ instead of $\Prw_{x,p_o}$. For transition probabilities
$p_\omega$ defined in terms of an environment $\omega$, we use the notation
$\Prw_{x,\omega}$. The corresponding expectation operators are denoted by
$\Erw_{x,p}$, $\Erw_x$ and $\Erw_{x,\omega}$, respectively. Every
$p\in\mathcal{P}$ gives in an obvious way rise to a homogeneous nearest
neighbor transition kernel on $\mathbb{Z}^d$, which we again denote by $p$.

For a subset $V\subset \mathbb{Z}^d$, we let $\tau_V = \inf\{n\geq 0 :
X_n\notin V\}$ be the first exit time from $V$, with $\inf \emptyset =
\infty$. 

Given $x,z\in\mathbb{Z}^d$, $p\in\mathcal{P}$ and a subset
$V\subset\mathbb{Z}^d$, we define
$$
\pi^{(p)}_V(x,z)= \Prw_{x,p}\left(X_{\tau_V}=z\right).
$$
For an environment $\omega\in\Omega$, we set
$$
\Pi_{V,\omega}(x,z)= \Prw_{x,\omega}\left(X_{\tau_V}=z\right).
$$
We mostly drop $\omega$ in the notation and interpret $\Pi_V(x,\cdot)$
as a random measure.

Recall the definitions of the sets $\mP$, $\mPsymone$ and $\mP_\e$ from the
introduction. For $0<\kappa<1/(2d)$, let
$$
\mPsym=\{p\in\mP_\kappa: p(e_i)=p(-e_i),\,i=1,\dots,d\},
$$
i.e. $\mPsym$ is the subset of $\mP_{\kappa}$ which contains all symmetric
probability distributions on $\{\pm e_i : i = 1,\dots,d\}$. The parameter
$\kappa$ was introduced in~\cite{BB} to bound the range of the
symmetric transition kernels under consideration. We can think of
$\kappa$ as a fixed but arbitrarily small number (the perturbation parameter
$\e$ is chosen afterward).

\subsubsection{Miscellaneous comments about notation}
Our constants are positive and depend only on the dimension $d\geq 3$
unless stated otherwise. In particular, they do {\it not} depend on
$L$, $p\in\mPsym$, $\omega$ or on any point $x\in\mathbb{Z}^d$.

By $C$ and $c$ we denote generic positive constants whose values can change
even in the same line. For constants whose values are fixed throughout a
proof we often use the symbols $K, C_1, c_1$. 

Many of our quantities, e.g. the transition kernels $\Ph_L$, $\ph_L$ or the
kernel $\Gamma_L$, are indexed by $L$. We normally drop the index in the
proofs. In contrast to~\cite{BB}, we do here not work with an additional
parameter $r$.

We often drop the superscript $(p)$ from notation and write $\pi_V$ for
$\pi^{(p)}_V$.  If $V=V_L$ is the ball around zero of radius $L$, we write
$\pi_L$ instead of $\pi_V$, $\Pi_L$ for $\Pi_V$ and $\tau_L$ for $\tau_{V}$.

By $\nP$ we denote sometimes a generic probability measure, and by
$\nE$ its corresponding expectation. If $A$ and $B$ are two events, we
often write $\nP(A;\, B)$ for $\nP(A\cap B)$.

If we write that a statement holds for ``$L$ large (enough)'', we
implicitly mean that there exists some 
$L_0>0$ depending only on the dimension
such that the statement is true for all $L\geq L_0$. This applies
also to phrases like ``$\delta$
(or $\e$, or $\kappa$) small (enough)''.

Some of our statements are only valid for large $L$ and $\e$ (or $\delta$,
or $\kappa$) sufficiently small, but we do not mention this every time.

\section{Results and concepts from the study of exit laws}
\label{SBB}
Our approach uses results and constructions from~\cite{BB}, where exit
measures from large balls under $\czero$ and $\cone(\e)$ are studied. We
adapt in this section those parts which will be frequently used in this
paper.  Some auxiliary statements from~\cite{BB}, which play here only a minor
role, will simply be cited when they are applied.

The overall idea of~\cite{BB} is to transport estimates on exit measures
inductively from one scale to the next, via a perturbation expansion for
the Green's function, which we recall first.

\subsection{A perturbation expansion}
Let $p =
\left(p(x,y)\right)_{x,y\in\mathbb{Z}^d}$ be a family of finite range
transition probabilities on $\mathbb{Z}^d$, and let $V\subset\mathbb{Z}^d$
be a finite set. The corresponding Green's kernel or Green's function for
$V$ is defined by
\begin{equation}
\label{prel-greensf}
g_V(p)(x,y) = \sum_{k=0}^\infty \left(1_Vp\right)^k(x,y).
\end{equation}
Now write $g$ for $g_V(p)$ and let $P$ be another transition
kernel with corresponding Green's function $G$ for $V$. With $\Delta=
1_V\left(P-p\right)$, the resolvent equation gives  
\begin{equation}
\label{prel-pbe1}
G -g = g\Delta G = G\Delta g.
\end{equation}
An iteration of~\eqref{prel-pbe1} leads to further expansions. Namely,
first one has
\begin{equation}
\label{prel-pbe2}
G -g = \sum_{k=1}^\infty \left(g\Delta\right)^kg, 
\end{equation}
and then, with $R= \sum_{k=1}^\infty\Delta^kp$, one arrives at
\begin{equation}
\label{prel-pbe3}
G = g\sum_{m=0}^\infty{\left(Rg\right)}^m\sum_{k=0}^\infty\Delta^k.
\end{equation}
We refer to~\cite{BB} for more details.

\subsection{Coarse grained transition kernels}
\label{SUBSprel-cgrw}
We fix once for all a probability density $\varphi\in
C^\infty(\mathbb{R}_+,\mathbb{R}_+)$ with compact support in $(1,2)$. Given
a transition kernel $p\in\mathcal{P}$ and a strictly positive function
$\psi = (m_x)_{x\in W}$, where $W\subset\mathbb{R}^d$, we define coarse
grained transition kernels on $W\cap\mathbb{Z}^d$ associated to
$(\psi,\,p)$,
\begin{equation}
\label{prel-cgpsi}
  \php_{\psi}(x,\cdot)  =
  \frac{1}{m_x}\int_{\mathbb{R}_+}
  \varphi\left(\frac{t}{m_x}\right)\pi_{V_t(x)}^{(p)}(x,\cdot)\dt
  t,\quad x\in W\cap\mathbb{Z}^d.
\end{equation}
Often $\psi\equiv m>0$ will be a constant, and then~\eqref{prel-cgpsi}
makes sense for all $x\in\mathbb{Z}^ d$ and therefore gives coarse grained
transition kernels on the whole grid $\mathbb{Z}^d$.

We now introduce particular coarse grained transition kernels for walking
inside the ball $V_L$, for both symmetric random walk and RWRE. We set up a
coarse graining scheme which will make the link between different scales
and allows us to transport estimates on mean sojourn times from one level
to the next. Our scheme is similar to that in~\cite{BB}, but does not
depend on an additional parameter $r$.

Let 
$$
s_L = \frac{L}{(\log L)^3} \quad\mbox{and}\quad r_L= \frac{L}{(\log L)^{15}}.
$$
We fix a smooth function $h :
\mathbb{R}_+\rightarrow\mathbb{R}_+$ satisfying
$$
  h(x)= \left\{\begin{array}{l@{\quad \mbox{for\ }}l}
      x & x\leq 1/2\\
      1 & x\geq 2\end{array}\right.,
$$
such that $h$ is concave and increasing on $(1/2,2)$.
Define $h_L : \overline{C}_L \rightarrow\mathbb{R}_+$ by
\begin{equation}
\label{eq-hlr}
 h_L(x) = \frac{1}{20}\max\left\{s_L
  h\left(\frac{\dL(x)}{s_L}\right),\, r_L\right\}.
\end{equation}
Note that in the setting of~\cite{BB}, this means that we always work with
the choice $r=r_L$, and there is no need keep $r$ in the notation.

We write $\Ph_L=\Ph_{L,\omega}$ for the coarse grained RWRE transition kernel inside
$V_L$ associated to $(\psi= \left(h_L(x)\right)_{x\in V_L},p_\omega)$,
\begin{equation}
\label{eq-cgkernel}
\Ph_L(x,\cdot) =
\frac{1}{h_L(x)}\int_{\mathbb{R}_+}
  \varphi\left(\frac{t}{h_L(x)}\right)\Pi_{V_t(x)\cap V_L}(x,\cdot)\dt
  t,
\end{equation}
and $\php_L$ for the corresponding kernel coming from symmetric random walk
with transition kernel $p\in\mathcal{P}$, where in the
definition~\eqref{eq-cgkernel} the random RWRE exit measure $\Pi$ is
replaced by $\pip$. For points $x\in\mathbb{Z}^d\backslash V_L$, we set
$\Ph_L(x,\cdot) =\php_L(x,\cdot) = \delta_x(\cdot)$.

Note our small abuse of notation: $\php_L$ is always defined in this way
and does {\it never} denote the coarse grained kernel~\eqref{prel-cgpsi}
associated to the constant function $\psi\equiv L$. Also note that $\Ph_L$
was denoted $\Ph_{L,r_L}$ in~\cite{BB}, and similarly $\ph_L$ was denoted
$\ph_{L,r_L}$. The kernel $\Ph_L$ is a random transition kernel depending
on $\omega$. However, when we consider $\Ph_L$ under $\Prw_{x,\omega}$,
then $\omega$ is fixed, but even in this case we usually write $\Ph_L$
instead of $\Ph_{L,\omega}$.

Two Green's function will play a crucial role (cf.~\eqref{prel-greensf}).
\begin{itemize}
\item $\Gh_L$ denotes the (coarse grained)
RWRE Green's function corresponding to $\Ph_L$.
\item $\ghp_L$ denotes the Green's function corresponding to $\php_L$.
\end{itemize}
The ``goodified'' version $\Ghg_L$ of $\Gh_L$ will be
introduced in Section~\ref{SUBSgoodifiedgreen}.

\subsection{Propagation of Condition \bf{C1}}

\label{SUBScspace}
We recapitulate in this part the technical Condition
$\cspace(\delta,L_0,L_1)$, which is propagated from one level to the next in~\cite{BB}.

\subsubsection{Assignment of transition kernels}
Let $L_0>0$ ($L_0$ shall play the role of a large constant). Define 
$L$-dependent symmetric transition kernels by
\begin{equation}
\label{kernelpL}
p_L(\pm e_i) =\left\{\begin{array}{l@{\quad \mbox{for\ }}l}
    1/(2d)& 0<L\leq L_0\\
    \frac{1}{2}\sum_{y\in\mathbb{Z}^d}\pE\left[\Ph_L(0,y)\right]\frac{y_i^2}{|y|^2}& L>L_0\end{array}\right..
\end{equation}

To be in position to formulate Condition $\cspace$, we recall some notation
from~\cite{BB}. Let
$\mathcal{M}_t$ be the set of smooth functions
$\psi:\mathbb{R}^d\rightarrow\mathbb{R}_+$ with first four derivatives
bounded uniformly by $10$ and
$$\psi\left(\left\{x\in\mathbb{R}^d: t/2<|x|<2t\right\}\right)\subset
(t/10, 5t).$$ For $p,q\in\mathcal{P}$ and $\psi\in\mathcal{M}_t$, define
\begin{align*}
D_{t,p,\psi,q}^{\ast} &=\sup_{x\in
  V_{t/5}}\left\|\left(\Pi_{V_t}-\pip_{V_t}\right)\ph_{\psi}^{(q)}(x,\cdot)\right\|_1,\\
 D_{t,p}^\ast &= \sup_{x\in V_{t/5}}\left\|\left(\Pi_{V_t}-\pip_{V_t}\right)(x,\cdot)\right\|_1.
\end{align*}
With $\delta > 0$, define for $i=1,2,3$ 
$$
b_i(L,p,\psi,q,\delta)= \pP\left(\left\{(\log L)^{-9+9(i-1)/4} < D_{L,p,\psi,q}^\ast \leq
  (\log L)^{-9 +9i/4}\right\}\cap\left\{D_{L,p}^\ast\leq \delta\right\}\right),
$$
and
$$
b_4(L,p,\psi,q,\delta)= \pP\left(\left\{D_{L,p,\psi,q}^\ast > (\log
    L)^{-3+3/4}\right\}\cup\left\{D_{L,p}^\ast> \delta\right\}\right).
$$
Let $\iota= (\log L_0)^{-7}$. Then Condition $\cspace$ is given as follows.
\subsubsection{Condition {\bf C1}}
Let $\delta > 0$ and $L_1\geq L_0\geq 3$. We say that $\cspace(\delta,L_0,L_1)$
holds if 
\begin{itemize}
\item  For all $3\leq L\leq 2L_0$, all $\psi\in\mathcal{M}_L$ and all $q\in\mPsymi$, 
$$
\pP\left(\left\{D_{L,p_o,\psi,q}^\ast > (\log L)^{-9}\right\}\cup\left\{D_{L,p_o}^\ast >
\delta\right\}\right) \leq \exp\left(-(\log (2L_0))^2\right).
$$
\item For all $L_0< L\leq L_1$, $L'\in [L,2L]$, $\psi\in\mathcal{M}_{L'}$ and 
  $q\in\mPsymi$, 
$$
b_i(L',p_L,\psi,q,\delta)\leq\frac{1}{4}\exp\left(-\left((3+i)/4\right)(\log
  L')^2\right)\quad\mbox{for }i=1,2,3,4.
$$
\end{itemize}

\subsubsection{The main statement for $\cspace$}

The first part of~\cite[Proposition 1.1]{BB} implies the following statement.
\begin{proposition}
  \label{main-prop-exitmeas}
  Assume $\cone$. For $\delta >0$ small enough, there exist
  $L_0=L_0(\delta)$ large and 
 $\e_0 = \e_0(\delta) > 0$ small with the following property: If $\e \leq \e_0$
  and $\czero(\e)$ is satisfied, then $\cspace\left(\delta,L_0,L\right)$ holds for
  all $L\geq L_0$.
\end{proposition}
For us, the important implication is that if $\cspace(\delta,L_0,L_1)$ is
satisfied, then for any $3\leq L\leq L_1$ and for all $L'\in[L,2L]$, all
$\psi\in\mathcal{M}_{L'}$ and all $q\in\mPsymi$,
\begin{equation}
\label{cspace-essential}
\pP\left(\left\{D_{L',p_L,\psi,q}^\ast > (\log L')^{-9}\right\}\cup\left\{D_{L',p_L}^\ast >
\delta\right\}\right) \leq \exp\left(-(\log L')^2\right).
\end{equation}
In~\cite[Lemma 2.2]{BB} it is shown that the transition kernels $p_L$ form
a Cauchy sequence. Their limit is given by the kernel $p_\infty$ defined
in~\eqref{eq-pinfty}, i.e. $\lim_{L\rightarrow}p_L(e_i)=p_{\infty}(e_i)$
for $i=1,\dots,d$. From this fact and the last display one can deduce that
the difference in total variation of the exit laws $\Pi_L$ and
$\pi_L^{(p_\infty)}$ is small when $L$ is large, in both a smoothed and
non-smoothed way. See Theorems 0.1 and 0.2 of~\cite{BB} for precise
statements. For us, it will be sufficient to keep in
mind~\eqref{cspace-essential} and the fact that
$\lim_{L\rightarrow}p_L=p_{\infty}$.

We follow the convention of~\cite{BB} and write ``assume
$\cspace(\delta,L_0,L_1)$'', if we assume $\cspace(\delta,L_0,L_1)$ for
some $\delta>0$ and some $L_1\geq L_0$, where $\delta$ can be chosen
arbitrarily small and $L_0$ arbitrarily large.

\subsection{Good and bad points}
\label{SUBSgoodandbad}
In~\cite[Section 2.2]{BB}, the concept of good and bad points inside
$V_L$ is introduced. It turns out that for controlling mean sojourn times,
we need a stronger notion of ``goodness'', see
Section~\ref{SUBSstgoodandbad}. It is however more convenient to first
recall the original classification.

Recall the assignment~\eqref{kernelpL}. A point $x\in V_L$ is {\it
  good} (with respect to $L$ and $\delta >0$), if
\begin{itemize}
\item For all $t\in[h_L(x),2h_L(x)]$, with $q=p_{h_L(x)}$,
  $$\left\|\left(\Pi_{V_t(x)}-\piq_{V_t(x)}\right)(x,\cdot)\right\|_1 \leq \delta.$$
\item If $\dL(x) > 2r$, then additionally
$$
\left\|\left(\Ph_L - \phq_L\right)\phq_L(x,\cdot)\right\|_1 \leq \left(\log
     h_L(x)\right)^{-9}.
$$
\end{itemize}
The set of environments where all points $x\in V_L$ are good is denoted 
$\good$.  A point $x\in V_L$ which is not good is called {\it bad}, and
the set of all bad points inside $V_L$ is denoted by
$\badP_L=\badP_L(\omega)$. 

\subsection{Goodified transition kernels and Green's function}
\label{SUBSgoodifiedgreen}
By replacing the coarse grained RWRE transition kernel at bad points $x\in
V_L$ by that of a homogeneous symmetric random walk, we obtain what we call
``goodified'' transition kernels inside $V_L$.

More specifically, write $p$ for $p_{s_L/20}$ stemming from the
assignment~\eqref{kernelpL}. The goodified
transition kernels are then defined as follows.
\begin{equation}
  \label{eq-goodifiedkernel}
  \Phg_L(x,\cdot) = \left\{\begin{array}{l@{\ \mbox{for\ }}l}
      \Ph_L(x,\cdot) & x \in V_L\backslash \badP_L\\
      \php_L(x,\cdot) & x \in \badP_L\end{array}\right..
\end{equation}
We write $\Ghg_L$ for the corresponding (random) Green's
function (denoted $\Ghg_{L,r_L}$ in~\cite{BB}).

Proposition~\ref{main-prop-exitmeas} will allow us to concentrate on
environments $\omega\in\good$, where $\Ph_L=\Phg_L$ and therefore also
$\Gh_L=\Ghg_L$. In the next section, we provide the necessary estimates for
the ``goodified'' coarse grained RWRE Green's function $\Ghg_L$.

\section{Control on Green's functions}
\label{Scontrolgreen}
We first recapitulate estimates on Green's functions
from~\cite{BB}. Then we
establish difference estimates for these functions, which will be used to
control differences of (quenched) mean sojourn times from balls $V_t(x)$ and
$V_t(y)$ that have a sufficiently large intersection.

\subsection{Bounds on Green's functions}
Recall that $\mPsym$ denotes the set of kernels which are symmetric in
every coordinate direction and $\kappa$-perturbations of the simple random
walk kernel. The statements in this section are valid for
small $\kappa$, in the sense that there exists $0<\kappa_0<1/(2d)$ such
that for $0<\kappa\leq \kappa_0$, the statements hold true for
$p\in\mPsym$, with constants that are uniform in $p\in\mPsym$.

Let $p\in\mPsym$ and $m\geq 1$. Denote by $\ph_{\psi_m}={\php}_{\psi_m}$
the coarse grained transition probabilities on $\mathbb{Z}^d$ associated to
$\psi_m = {(m_x)}_{x\in\mathbb{Z}^d}$, where $m_x= m$ is chosen constant in
$x$, cf.~\eqref{prel-cgpsi}. We mostly drop $p$ from notation. The
kernel $\ph_{\psi_m}$ is centered, with covariances
$$
\sum_{y\in\mathbb{Z}^d}(y_i-x_i)(y_j-x_j)\ph_{\psi_m}(x,y) = \lambda_{m,i}\delta_{i}(j),
$$
where for large $m$, 
$C^{-1}<\lambda_{m,i}/m^2< C$ for some $C>0$.
We set 
$$
{\bf \Lambda}_m= \left(\lambda_{m,i}\delta_i(j)\right)_{i,j=1}^d,\quad \mJ(x) =
|{\bf \Lambda}_m^{-1/2}x|\quad\hbox{for } x\in\mathbb{Z}^d,
$$
and denote by  
$$\gh_{m,{\mathbb{Z}^d}}(x,y) =
\sum_{n=0}^\infty(\ph_{\psi_m})^n(x,y)$$ 
the Green's function corresponding
to $\ph_{\psi_m}$. In~\cite{BB}, the following behavior of 
$\gh_{m,{\mathbb{Z}^d}}$ was established. The proof is based on a local central limit theorem for $\ph_{\psi_m}$,
which we do not restate here.
\begin{proposition}
\label{super-behaviorgreen}
Let $p\in\mPsym$. Let $x,y\in\mathbb{Z}^d$, and assume $m \geq m_0 > 0$ large enough. 
\begin{enumerate}
\item For $|x-y| < 3m$,
$$
\gh_{m,{\mathbb{Z}^d}}(x,y) = \delta_{x}(y) + O(m^{-d}).
$$
\item For $|x-y| \geq 3m$, there exists a constant $c(d) > 0$ such that
$$
  \gh_{m,{\mathbb{Z}^d}}(x,y)= \frac{c(d)\det{\bf \Lambda}_m^{-1/2}}{\mJ(x-y)^{d-2}} +
  O\left(\frac{1}{|x-y|^{d}}\left(\log\frac{|x-y|}{m}\right)^d\right). 
$$
\end{enumerate}
\end{proposition}
Recall that $\tau_L=\tau_{V_L}$ denotes the first exit time from $V_L$.
The proposition can be used to estimate the corresponding Green's function
for $V_L$,
$$\gh_{m, V_L}(x,y) =
\sum_{n=0}^\infty\left(1_{V_L}\ph_{\psi_m}\right)^n(x,y).$$ Indeed,
$\gh_{m, V_L}$ is bounded from above by $\gh_{m,{\mathbb{Z}^d}}$, and more
precisely, the strong Markov property shows
\begin{equation}
\label{localclt-greenonball}
  \gh_{m,V_L}(x,y) = \Erw_{x,\ph_{\psi_m}}\left[\sum_{k=0}^{\tau_L-1}1_{\{X_k = y\}}\right] =
  \gh_{m,\mathbb{Z}^d}(x,y) - \Erw_{x,\ph_{\psi_m}}\left[\gh_{m,\mathbb{Z}^d}\left(X_{\tau_L},y\right)\right]. 
\end{equation}
Here, according to our notational convention, $\Erw_{x,\ph_{\psi_m}}$
is the expectation with respect to $\Prw_{x,\ph_{\psi_m}}$, the law  
of a random walk started at $x$ and running with kernel $\ph_{\psi_m}$.

We next recall the definition of the (deterministic) kernel $\Gamma_L$, which was
introduced in~\cite{BB} to dominate coarse grained Green's
functions from above. 

Since we always work with $r=r_L$, we write $\Gamma_L$ instead of
$\Gamma_{L,r}$ as in~\cite{BB}, and in the proofs, the index $L$ is dropped
as well. We formulate our definitions and results in terms of the larger
ball $V_{L+r_L}$, so that we can refer to the proofs in~\cite{BB}. For
$x\in V_{L+r_L}$, let
$$
  \dtL(x) = \max\left(\frac{\dist_{L+r_L}(x)}{2},3r_L\right),\quad a(x) =\min\left(\dtL(x),s_L\right).
$$
For $x,y\in V_{L+r_L}$, the kernel $\Gamma_L$ is now defined by 
\begin{equation}
\label{defgamma}
  \Gamma_L(x,y) = \min\left\{\frac{\tilde{d}(x)\tilde{d}(y)}{a(y)^2(a(y) +
    |x-y|)^d},\,\frac{1}{a(y)^2(a(y) + |x-y|)^{d-2}}\right\}.
\end{equation}
For $x\in V_{L+r_L}$, we write $U(x) =
V_{a(x)}(x)\cap V_{L+r_L}$ for the $a(x)$-neighborhood around $x$. Given two positive functions $F, G:
V_{L+r_L}\times V_{L+r_L}\rightarrow\mathbb{R}_+$, we write $F\preceq G$ if
for all $x,y\in V_{L+r_L}$,
$$
  F(x,U(y)) \leq G(x,U(y)),
$$
where $F(x,U)$ stands for $\sum_{y\in U\cap\mathbb{Z}^d}F(x,y)$. We write
$F\asymp 1$, if there is a constant $C>0$ such that for all $x,y\in
V_{L+r_L}$,
$$
\frac{1}{C}F(x,y)\leq F(\cdot,\cdot)\leq CF(x,y)\quad\mbox{on }U(x)\times U(y).
$$
We shall repeatedly need some properties of $\Gamma_L$, which form part
of~\cite[Lemma 4.4]{BB}.

\begin{lemma}[Properties of $\Gamma_L$]\
  \label{super-gammalemma}
  \begin{enumerate}
  \item $\Gamma_L\asymp 1$.
  \item For $1\leq j \leq \frac{1}{3r}s_L$, with $\mathcal{E}_j=\{y\in V_{L+r_L} : \dtL(y)\leq 3jr_L\}$,
$$
  \sup_{x\in V_{L+r_L}}\Gamma_L(x,\mathcal{E}_j) \leq C\log(j+1),
$$
and for $0\leq \alpha < 3$,
$$
  \sup_{x\in V_{L+r_L}}\Gamma_L\left(x,\sh_L\left(s_L, L/(\log
      L)^\alpha\right)\right) \leq C(\log\log L)(\log L)^{6-2\alpha}. 
$$
\item For $x\in V_{L+r}$,
$$
  \Gamma_L(x,V_L) \leq C\max\left\{\frac{\dtL(x)}{L}(\log L)^6,\,
    \left(\frac{\dtL(x)}{r_L}\wedge \log\log L\right)\right\}.
$$
\end{enumerate}
\end{lemma}
We now formulate the key estimate, which shows how both $\ghq_L$
and $\Ghg_L$ can be dominated from above by the deterministic kernel $\Gamma_L$.
See~\cite[Lemma 4.2]{BB} for a proof.
\begin{lemma}\ 
\label{superlemma} 
\begin{enumerate}
\item There exists a constant $C>0$ such that for all $q\in\mPsym$,
$$
      \ghq_L\preceq C\Gamma_L.
    $$
  \item Assume $\cspace(\delta,L_0,L_1)$, and let $L_1\leq L\leq L_1(\log
    L_1)^2$. There exists a constant $C>0$ such that for $\delta>0$ small,
$$
  \Ghg_L \preceq C\Gamma_L.
$$  
\end{enumerate}
\end{lemma}

\subsection{Difference estimates}
\label{super-difference}
For controlling mean sojourn times, we will need difference estimates for
the coarse grained Green's functions $\ghq_L$ and $\Ghg_L$. We first recall
our notation:
\begin{itemize}
\item For $q\in\mPsym$, $\ghq_L$ is the Green's function in $V_L$ corresponding to
$\phq_L$.
\item $\Ghg_L$ is the Green's function in $V_L$ corresponding to
  $\Phg_L$, cf.~\eqref{eq-goodifiedkernel}.
 
\item For $m>0$, $\ghq_{m,V_L}$ is the Green's function in $V_L$ corresponding
  to $1_{V_L}\phq_{\psi_m}$, where $\psi_m\equiv m$.
  \item $\ghq_{m,\mathbb{Z}^d}$ is the Green's function on
    $\mathbb{Z}^d$ corresponding to $\phq_{\psi_m}$, where $\psi_m\equiv m$.
  \end{itemize}

\begin{lemma}
\label{super-greendifference}
There exists a constant $C >0$ such that for all $q\in\mPsym$,
\begin{enumerate}
\item
$$
\sup_{x,x'\in V_L: |x-x'|\leq s_L}\sum_{y\in V_L}\big|\ghq_L(x,y)-\ghq_L(x',y)\big|
\leq C(\log\log L)(\log L)^3.
$$
\item Assume $\cspace(\delta,L_0,L_1)$, and let $L_1\leq L\leq L_1(\log
    L_1)^2$. There exists a constant $C>0$ such that for $\delta>0$ small,
$$
\sup_{x,x'\in V_L: |x-x'|\leq s_L}\sum_{y\in V_L}\big|\Ghg_L(x,y)-\Ghg_L(x',y)\big|
\leq C(\log\log L)(\log L)^3.
$$
\end{enumerate}
\end{lemma}
\begin{proof1}
  (i) The underlying one-step transition kernel is always given by
  $q\in\mPsym$, which we constantly omit in this proof, i.e. $\ph_{\psi_m}
  = \phq_{\psi_m}$, $\gh_{m,V_L} = \ghq_{m,V_L}$, $\gh_{m,\mathbb{Z}^d} =
  \ghq_{m,\mathbb{Z}^d}$, or $\Prw_{x} = \Prw_{x,q}$, and so on. Also, we
  suppress the subscript $L$, i.e. $\gh=\gh_L$. Set $m = s_L/20$. We write
\begin{align}
\label{super-greendifference-eq0}
  \lefteqn{\sum_{y\in V_L}\left|\gh(x,y)-\gh(x',y)\right|}\\
  &\leq \sum_{y\in
      V_L}\left|\left(\gh-\gh_{m,V_L}\right)(x,y)\right| + \sum_{y\in
      V_L}\left|\gh_{m,V_L}(x,y)-\gh_{m,V_L}(x',y)\right| + \sum_{y\in
      V_L}\left|\left(\gh_{m,V_L}-\gh\right)(x',y)\right|.\nonumber
\end{align}
If $x\in V_L\backslash \sh_L(2s_L)$, we have $\ph(x,\cdot)=
\ph_{\psi_m}(x,\cdot)$. Clearly, $\sup_{x\in V_L}\gh_{m,V_L}(x,\sh_L(2s_L))
\leq C$. Thus, with $\Delta= 1_{V_L}\left(\ph_{\psi_m}-\ph\right)$, the
perturbation expansion~\eqref{prel-pbe1} and
Lemma~\ref{super-gammalemma} yield (remember $\gh\preceq C\Gamma$ by
Lemma~\ref{superlemma})
\begin{align*} \sum_{y\in V_L}\left|(\gh_{m,V_L}-\gh)(x,y)\right|&=
  \sum_{y\in V_L}|\gh_{m,V_L}\Delta\gh(x,y)|\\ 
&\leq 2\,\gh_{m,V_L}(x,\sh_L(2s_L))\sup_{v\in
    \sh_L(3s_L)}\gh(v,V_L) \leq C(\log L)^3.
\end{align*}
It remains to handle the middle term
of~\eqref{super-greendifference-eq0}. By~\eqref{localclt-greenonball}, 
\begin{align*}
  \lefteqn{\gh_{m,V_L}(x,y)-\gh_{m,V_L}(x',y)}\\
  & = \gh_{m,\mathbb{Z}^d}(x,y)-
  \gh_{m,\mathbb{Z}^d}(x',y) +
  \Erw_{x',\ph_m}\left[\gh_{m,\mathbb{Z}^d}(X_{\tau_L},y)\right] -
  \Erw_{x,\ph_m}\left[\gh_{m,\mathbb{Z}^d}(X_{\tau_L},y)\right].  
\end{align*}
Using Proposition~\ref{super-behaviorgreen}, it follows that for $|x-x'|
\leq s_L$,
$$
  \sum_{y\in
      V_L}\left|\gh_{m,\mathbb{Z}^d}(x,y)-\gh_{m,\mathbb{Z}^d}(x',y)\right| \leq
  C(\log L)^3. 
$$
At last, we claim that
\begin{equation}
\label{super-greendifference-eq1}
  \sum_{y\in V_L}\left|\Erw_{x',\ph_{\psi_m}}\left[\gh_{m,\mathbb{Z}^d}(X_{\tau_L},y)\right] -
    \Erw_{x,\ph_{\psi_m}}\left[\gh_{m,\mathbb{Z}^d}(X_{\tau_L},y)\right]\right|
  \leq C(\log\log L){(\log L)}^3.  
\end{equation}
Since $|x-x'| \leq m$, we can define on the same probability space, whose
probability measure we denote by $\mathbb{Q}$, a random walk $(Y_n)_{n\geq
  0}$ starting from $x$ and a random walk $(\tilde{Y}_n)_{n\geq 0}$
starting from $x'$, both moving according to $\ph_{\psi_m}$ on
$\mathbb{Z}^d$, such that for all times $n$, $|Y_n - \tilde{Y}_n| \leq
s_L$. However, with $\tau= \inf\{n \geq 0 : Y_n \notin V_L\}$,
$\tilde{\tau}$ the same for $\tilde{Y}_n$, we cannot deduce that
$|Y_{\tau}-\tilde{Y}_{\tilde{\tau}} | \leq s_L$, since it is
possible that one of the walks, say $Y_n$, exits $V_L$ first and then moves far
away from the exit point, while the other walk  
$\tilde{Y}_n$ might still be inside $V_L$. In order to show that
such an event has a small probability, we argue similarly
to~\cite[Proposition 7.7.1]{LawLim}. Define
$$
\sigma(s_L) = \inf\left\{n\geq 0 : Y_n\in \sh_L(s_L)\right\},
$$
and analogously $\tilde{\sigma}(s_L)$. Let $\vartheta =
\sigma(s_L)\wedge\tilde{\sigma}(s_L)$.  Since
$|Y_{\vartheta}-\tilde{Y}_{\vartheta}| \leq s_L$, 
$$
\max\left\{\sigma(2s_L),\,\tilde{\sigma}(2s_L)\right\} \leq \vartheta.
$$
For $k\geq 1$, we introduce the events
\begin{align*}
B_k&=\left\{\big|Y_i-Y_{\sigma(2s_L)}\big| > ks_L\mbox{ for
      all } i=\sigma(2s_L),\dots,\tau\right\},\\
\tilde{B}_k&=\left\{\big|\tilde{Y}_i-\tilde{Y}_{\tilde{\sigma}(2s_L)}\big| > ks_L\mbox{ for
      all } i=\tilde{\sigma}(2s_L),\dots,\tilde{\tau}\right\}.
\end{align*}
By the strong Markov property and the gambler's ruin estimate of~\cite{LawLim},
p. 223 (7.26), there exists a constant $C_1>0$ independent of $k$ such that 
$$
\mathbb{Q}\left(B_k\cup \tilde{B}_k\right) \leq C_1/k
$$
for some $C_1>0$ independent of $k$. Applying the triangle inequality to
$$
  Y_{\tau} - \tilde{Y}_{\tilde{\tau}} =
  (Y_{\tau}-Y_{\vartheta}) +
  (Y_{\vartheta}-\tilde{Y}_{\vartheta}) + (\tilde{Y}_{\vartheta}-\tilde{Y}_{\tilde{\tau}}),
$$
we deduce, for $k \geq 3$,  
$$
\mathbb{Q}\left(\big|Y_\tau - \tilde{Y}_{\tilde{\tau}}\big|\geq
  ks_L\right) \leq 2C_1/(k-1).
$$
Since $|Y_\tau-\tilde{Y}_{\tilde{\tau}}| \leq 2(L +s_L) \leq 3L $, it follows that
$$
\mathbb{E}_{\mathbb{Q}}\left[\big|Y_\tau-\tilde{Y}_{\tilde{\tau}}\big|\right] \leq
\sum_{k=1}^{3L}\mathbb{Q}\left(\big|Y_\tau-\tilde{Y}_{\tilde{\tau}}\big|
  \geq k\right)\leq C(\log\log L)s_L.
$$
Also, for $v,w$ outside and $y$ inside $V_L$,
$$
  \left|\frac{1}{|v-y|^{d-2}}-\frac{1}{|w-y|^{d-2}}\right|
  \leq 
  C\frac{|v-w|}{(L+1-|y|)^{d-1}}.  
$$
Applying
Proposition~\ref{super-behaviorgreen},~\eqref{super-greendifference-eq1}
now follows from the last two displays and a summation over
$y\in V_L$.\\
(ii) We take $p=p_{s_L/20}$ stemming from the
assignment~\eqref{kernelpL} and work with $p$ as the underlying one-step
transition kernel, which will be suppressed from the notation,
i.e. $\ph=\php$ and $\gh=\ghp$.

Let $x, x'\in V_L$ with $|x-x'|\leq s_L$ and set $\Delta =
1_{V_L}(\Phg-\ph)$. With $B = V_L\backslash\sh_L(2r_L)$,  
$$ 
\Ghg = \gh 1_{B}\Delta\Ghg +\gh 1_{V_L\backslash B}\Delta\Ghg +\gh.
$$
Replacing successively $\Ghg$ in the first summand on the right-hand
side,
$$
\Ghg= \sum_{k=0}^\infty{\left(\gh1_B\Delta\right)}^k\gh +
    \sum_{k=0}^\infty{\left(\gh1_B\Delta\right)}^k\gh1_{V_L\backslash B}\Delta\Ghg
    = F+F1_{V_L\backslash B}\Delta\Ghg,  
$$
where we have set $F=\sum_{k=0}^\infty{\left(\gh1_B\Delta\right)}^k\gh$.
With $R = \sum_{k=1}^{\infty}(1_B\Delta)^k\ph$,
expansion~\eqref{prel-pbe3} gives
\begin{equation}
\label{super-diffF1}
F = \gh\sum_{m=0}^\infty
(R\gh)^m\sum_{k=0}^\infty\left(1_B\Delta\right)^k =
\gh\sum_{k=0}^\infty\left(1_B\Delta\right)^k + \gh RF.
\end{equation}
By the arguments given in the proof of Lemma~\ref{superlemma} (ii)
in~\cite{BB} (note $\|1_B\Delta\|_1\leq \delta$, and $\|1_B\Delta\ph\|_1\leq C(\log
L)^{-9}$), one deduces $|F|\preceq C\Gamma$. By
Lemma~\ref{super-gammalemma} (ii) and (iii), we then see that for large $L$,
uniformly in $x\in V_L$,
$$
\|F1_{V_L\backslash B}\Delta\Ghg(x,\cdot)\|_1 \leq
C\Gamma(x,\sh_L(2r_L))\sup_{v\in\sh_L(3r_L)}\Gamma(v,V_L) \leq C\log\log L.
$$
Therefore,
$$
  \sum_{y\in V_L}\big|\Ghg(x,y)-\Ghg(x',y)\big| \leq C\log\log L +
  \sum_{y\in V_L}\left|F(x,y)-F(x',y)\right|.
$$
Using~\eqref{super-diffF1} and twice part (i),
\begin{align}
\label{super-diffF2}
\lefteqn{\sum_{y\in V_L}\left|F(x,y)-F(x',y)\right|}\nonumber\\ 
&\leq & \sum_{y\in
  V_L}\Big|\gh\sum_{k=0}^\infty\left(1_B\Delta\right)^k(x,y)
  -\gh\sum_{k=0}^\infty\left(1_B\Delta\right)^k(x',y)\Big| + \sum_{y\in V_L}\left|\gh RF(x,y) - \gh RF(x',y)\right|.
\end{align}
The first expression on the right is estimated by
$$
\sum_{y\in V_L}\Big|\sum_{w\in
    V_L}\left(\gh(x,w)-\gh(x',w)\right)\sum_{k=0}^\infty\left(1_B\Delta\right)^k(w,y)\Big|
\leq C(\log\log L)(\log L)^3,
$$
where we have used part (i) and $\|1_B\Delta(w,\cdot)\|_1\leq \delta$.
The second factor of~\eqref{super-diffF2} is again bounded by (i) and the fact
that for $u \in V_L$,
\begin{align*}
\sum_{y\in V_L}|RF(u,y)| &= \sum_{y\in
  V_L}\Big|\sum_{k=1}^\infty\left(1_B\Delta\right)^k\ph F(u,y)\Big|\\
& \leq \sum_{k=0}^\infty \left\|1_B\Delta(u,\cdot)\right\|_1^k\sup_{v\in
  B}\left\|1_B\Delta\ph(v,\cdot)\right\|_1\sup_{w\in V_L}\sum_{y\in V_L}|F(w,y)|\\ 
&\leq C(\log L)^{-9+6} = C (\log L)^{-3}.
\end{align*}
Altogether, this proves part (ii).
\end{proof1}

\section{Mean sojourn times}
\label{Stimes}
Using the results about the exit measures from 
Proposition~\ref{main-prop-exitmeas} and the estimates for the Green's
functions from Section~\ref{Scontrolgreen}, we proof in this part our main
results on mean sojourn times in large balls.

\subsection{Condition {\bf C2} and the main technical statement}

\label{SUBSctimes}
Similarly to Condition $\cspace(\delta,L_0,L_1)$, cf. Section~\ref{SUBScspace},
we formulate a condition on the mean sojourn times which we propagate from one
level to the next.

We first introduce a monotone increasing function which will upper and
lower bound the
normalized mean sojourn time in the ball. Let $0<\eta<1$, and define
$f_\eta: \mathbb{R}_+\rightarrow \mathbb{R}_+$ by setting
$$
  f_\eta(L)= \frac{\eta}{3}\sum_{k=1}^{\lceil \log L\rceil}k^{-3/2}.
$$
Note $\eta/3\leq f_\eta(L)< \eta$ and therefore $\lim_{\eta\downarrow 0}\lim_{L\rightarrow\infty}f_\eta(L) = 0$.

Recall that $\Erw_x=\Erw_{x,p_o}$ is the expectation with respect to simple random walk
starting at $x\in\mathbb{Z}^d$, and $\tau_L$ is the first exit time from $V_L$.

\subsubsection{Condition {\bf C2}}
We say that $\ctime(\eta,L_1)$ holds, if for all
$3\leq L \leq L_1$,
$$
\pP\left(\Erw_{0,\omega}\left[\tau_L\right] \notin \left[1-f_\eta(L),\,
      1+f_\eta(L)\right]\cdot\Erw_0\left[\tau_L\right]\right) \leq
\exp\left(-(1/2)(\log L)^{-2}\right).
$$

Our main technical result for the mean sojourn times is
\begin{proposition}
\label{main-prop-times}
  Assume $\cone$ and $\cB$, and let $0<\eta<1$. There exists
  $\e_0=\e_0(\eta)> 0$ with the following
  property: If $\e \leq \e_0$ and $\czero(\e)$ holds, then there exists $L_0=L_0(\eta)>0$ such that for $L_1\geq L_0$,
$$
\ctime(\eta,L_1)\Rightarrow \ctime(\eta,L_1(\log L_1)^2).
$$
\end{proposition}
\begin{remark}
  Given $\eta$ and $L_0$, we can always guarantee (by making $\e$ smaller if
  necessary) that $\czero(\e)$ implies
  $\ctime(\eta,L_0)$.
 \end{remark} 

The proof of this statement is deferred to Section~\ref{SUBSctimes-proofs}.

 \subsection{Some preliminary results}
We begin with an elementary statement about the mean time a symmetric
random walk with kernel $p\in\mPsym$ spends in the ball $V_L$.
\begin{lemma}
\label{LmeanRWp}
Let $p\in\mPsym$, and let $x\in V_L$. Then
$$
L^2-|x|^2\leq \Erw_{x,p}\left[\tau_L\right]\leq (L+1)^2-|x|^2.
$$
\end{lemma}
The proof of this standard lemma (see e.g.~\cite[Proposition
6.2.6]{LawLim}) uses the fact that $|X_{n\wedge\tau_L}|^2-n\wedge \tau_L$ is
a martingale, which leads by optional stopping to $\Erw_{x,p}[\tau_L] =
\Erw_{x,p}[|X_{\tau_L}|^2] -|x|^2$.  In particular, for different
$p,q\in\mPsym$, the corresponding mean sojourn times satisfy
$$
\Erw_{0,p}\left[\tau_L\right] = \Erw_{0,q}\left[\tau_L\right](1+O(L^{-1})).
$$
We will compare the RWRE sojourn times on all scales with
$\Erw_{0}[\tau_L]$, the corresponding quantity for simple random walk. This
is somewhat in contrast to our comparison of the exit measure in~\cite{BB},
where we use the scale-dependent kernels $p_L$ given by~\eqref{kernelpL}.

Using that $\mu$ is supported on transition probabilities which are
balanced in the first coordinate direction, we obtain a similar upper
bound for the mean sojourn time of the RWRE.
\begin{lemma}
\label{LmeanRWREapriori}
For $\omega\in (\mathcal{P}_\e)^{\mathbb{Z}^d}\cap(\mPsymone)^{\mathbb{Z}^d}$,
$$
\Erw_{x,\omega}\left[\tau_L\right]\leq \frac{d}{1-2\e d}(L+1)^2-(x\cdot e_1)^2.
$$
\end{lemma}
\begin{proof1}
  For $\omega\in(\mPsymone)^{\mathbb{Z}^d}$, $\omega_x(e_1) =
  \omega_x(-e_1)$ for each $x\in\mathbb{Z}^d$. Then
$$
M_n = (X_n\cdot e_1)^2 - \sum_{k=0}^{n-1}\left(\omega_{X_k}(e_1) +
  \omega_{X_k}(-e_1)\right)
$$
is a $\Prw_{x,\omega}$-martingale with respect to the
filtration generated by the walk $(X_n)_{n\geq 0}$. By the optional stopping
theorem, $\Erw_{x,\omega}\left[M_{n\wedge\tau_L}\right] = (x\cdot e_1)^2$.
Since for $\omega\in(\mathcal{P}_\e)^{\mathbb{Z}^d}$,
$$
\omega_{X_k}(e_1)+\omega_{X_k}(-e_1) \geq 1/d -2\e,$$ 
it follows that
$$
\Erw_{x,\omega}\left[n\wedge\tau_L\right]\leq
{\left(1/d-2\e\right)}^{-1}\Erw_{x,\omega}\left[(X_{n\wedge\tau_L}\cdot
  e_1)^2\right] - (x\cdot e_1)^2.
$$
Letting $n\rightarrow\infty$ proves the statement.
\end{proof1}
\begin{remark}
  Conditions $\czero(\e)$ and $\cB$ guarantee that the event
  $(\mathcal{P}_\e)^{\mathbb{Z}^d}\cap(\mPsymone)^{\mathbb{Z}^d}$ has full
  $\pP$-measure. The a priori fact that $\Erw_{0,\omega}[\tau_L]\leq CL^2$
  for almost all environments will be crucial to obtain more precise
  bounds on these times. 
\end{remark}

We will now express the mean sojourn time of the RWRE in $V_L$ in terms of
mean sojourn times in smaller balls $V_t(x)\subset V_L$, for $t\in[h_L(x),2h_L(x)]$.
\begin{figure}[ht]
\begin{center}\parbox{5.5cm}{\includegraphics[width=5cm]{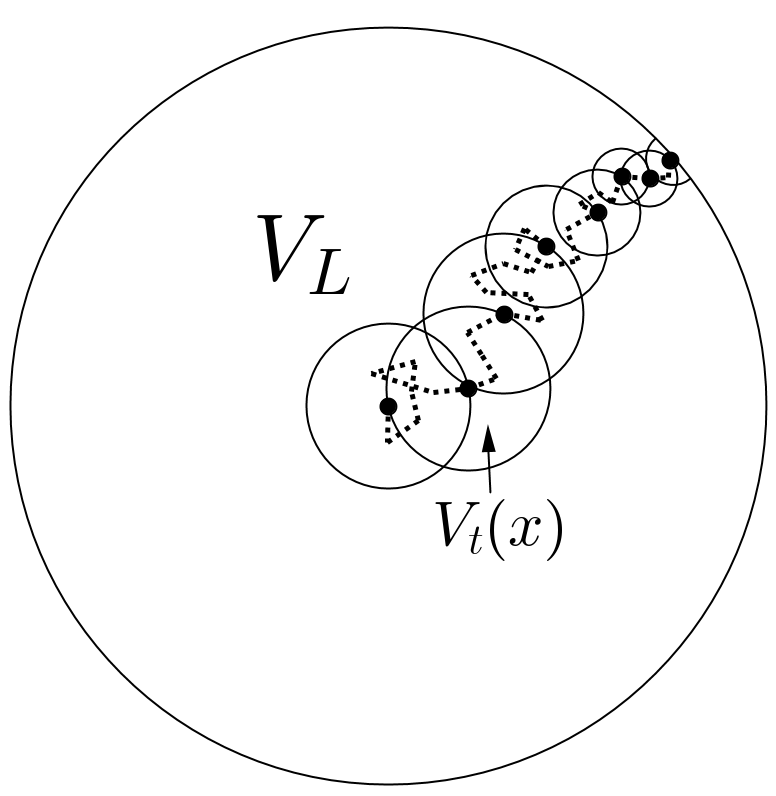}}
\parbox{10cm}{
  \caption{We express the mean sojourn time $\Erw_{0,\omega}[\tau_L]$ as a
    convolution of the coarse grained RWRE Green's function $\Gh_L$ with
    mean sojourn times in smaller balls $V_t(x)\cap V_L$, where
    $t\in[h_L(x),2h_L(x)]$ (see Lemma~\ref{times-keylemma}). Inductive control over the sojourn times on
    smaller scales $\leq s_L$ and over the Green's function then allow us
    to obtain the right estimate for $V_L$.}}
\end{center}
\end{figure}
Recall the definition of $h_L$ and the corresponding coarse graining scheme
inside $V_L$. As in Section~\ref{SUBSprel-cgrw}, we put
$$
s_t = \frac{t}{(\log t)^3} \quad\mbox{and}\quad r_t= \frac{t}{(\log t)^{15}}.
$$
Let $h_t^x(\cdot)= h_{t}(\cdot-x)$, where $h_{t}(\cdot-x)$ is defined
in~\eqref{eq-hlr} (with $L$ replaced by $t$). By translating the
origin into $x$, we transfer the coarse graining schemes on $V_L$ in the
obvious way to $V_t(x)$, using $h_t^x$ instead of $h_L$. We write $\Ph_t^x$
for the coarse grained transition probabilities in $V_t(x)$ associated to
$((h_t^x(y))_{y\in V_t(x)},p_\omega)$, cf.~\eqref{eq-cgkernel}. Given
$p\in\mPsym$, the kernel $\ph_t^{(p),x}$ is defined similarly, with
$p_\omega$ replaced by $p$.

For the corresponding Green's functions we use the expressions $\Gh_t^x$
and $\gh_t^{(p),x}$.  If we do not keep $x$ as an index, we always mean
$x=0$ as before. If it is clear with which $p$ we are working, we drop the
superscript $(p)$.  Notice that for $y,z \in V_t(x)$ and $p\in\mPsym$, we
have $\ph_t^{(p),x}(y,z) = \php_t(y-x,z-x)$ and $\ghp_t(y,z) =
\gh_t^{(p)}(y-x,z-x)$. Since $p_\omega$ is in general not homogeneous in
space, this is not true for $\Ph_t^x$ and $\Gh_t^x$.

Define the ``coarse grained'' RWRE sojourn time
$$
  \Lambda_L(x) = 1_{V_L}(x)\, 
  \frac{1}{h_L(x)}\int\limits_{\mathbb{R}_+}\varphi\left(\frac{t}{h_L(x)}\right)\Erw_{x,\omega}\left[\tau_{V_t(x)\cap
      V_L}\right]\dt t,    
$$
and the analog for random walk with kernel $p\in\mPsym$,
$$
  \lambdap_L(x) = 1_{V_L}(x)\,
      \frac{1}{h_L(x)}\int\limits_{\mathbb{R}_+}\varphi\left(\frac{t}{h_L(x)}\right)\Erw_{x,p}\left[\tau_{V_t(x)\cap
        V_L}\right]\dt t.  
$$ 
We also consider the corresponding quantities $\Lambda_t^x$ and $\lambda_t^{(p),x}$ for
balls $V_t(x)$. For example, 
$$
\Lambda_t^x(y) = 1_{V_t(x)}(y)\,\frac{1}{h_t^x(y)}\int\limits_{\mathbb{R}_+}\varphi\left(\frac{s}{h_t^x(y)}\right)\Erw_{y,\omega}\left[\tau_{V_s(y)\cap
  V_t(x)}\right]\dt s.  
$$
Note that we should rather write $\Lambda_{L,\omega}$ and
$\Lambda_{t,\omega}^x$, but we again suppress $\omega$ in the notation.
In the rest of this part, we often let operate kernels on mean sojourn times from the left. As an example,
$$
\Gh_L\Lambda_L(x) = \sum_{y\in V_L}\Gh_L(x,y)\Lambda_L(y).
$$
Both $\Lambda_L$ and $\Gh_L$ should be understood as random functions, but
sometimes (for example in the proof of the next statement) the environment
$\omega$ is fixed.

The basis for our inductive scheme is now established by 
\begin{lemma}
\label{times-keylemma}
For environments $\omega\in\left(\mathcal{P}_\varepsilon\right)^{\mathbb{Z}^d}$,
$x\in\mathbb{Z}^d$, 
$$\Erw_{x,\omega}\left[\tau_L\right] = \Gh_L\Lambda_L(x).$$
In particular, for $\omega$ the homogeneous environment with transition
probabilities given by $p\in\mPsym$, 
$$\Erw_{x,p}\left[\tau_L\right] = \ghp_L\lambdap_L(x).$$
\end{lemma}
\begin{proof1}
  We construct a probability space where we can observe in $V_L$ both the
  random walk running with transition kernel $p_\omega$ and its coarse
  grained version running with kernel $\Ph_L(\omega)$.  In this direction,
  we take a probability space $(\Xi, \mathcal{A}, \pQ)$ carrying a family
  of i.i.d. $[1,2]$-valued random variables $(\xi_n: n\in \mathbb{N})$
  distributed according to $\varphi(t)\dt t$.  We now consider the
  probability space $((\mathbb{Z}^d)^{\mathbb{N}}\times\Xi,
  \pG\otimes\mathcal{A}, \Prw_{x,\omega}\otimes\pQ)$. By a small abuse of
  notation, we denote here by $X_n$ the projection on the $n$th component
  of the first factor of $(\mathbb{Z}^d)^{\mathbb{N}}\times \Xi$, so that
  under $\Prw_{x,\omega}\otimes\pQ$, $(X_n)_{n\geq 0}$ evolves as the
  canonical Markov chain under $\Prw_{x,\omega}$.

  Set $T_0 = 0$ and define the ``randomized'' stopping times
$$
T_{n+1} = \inf\left\{m> T_n: X_m\notin V_{\xi_{T_n}\cdot
    h_L(X_{T_n})}\left(X_{T_n}\right)\right\}\wedge \tau_L.
$$
Then the coarse grained Markov chain in $V_L$ running with transition
kernel $\Ph_{L,omega}$ can be obtained by observing $X_n$ at the times
$T_n$, that is by considering $(X_{T_n})_{n\geq 0}$. Moreover, the Markov
property of $X_n$ and the i.i.d. property of the $\xi_n$ ensure that under
$\Prwt_{x,\omega}$, conditionally on $X_{T_n}$, the random vector
$((X_{T_n},X_{T_n+1},\dots),T_{n+1}-T_n)$ is distributed as
$((X_{0},X_{1},\dots),T_{1})$ under $\Prwt_{X_{T_n},\omega}$. Indeed,
formally one may define the filtration $\pG_n = \sigma(X_0,\dots,X_n,\xi_0, \dots,
\xi_{n-1})$. Then $(X_n)_{n\geq 0}$ is also a $\pG_n$-Markov chain. By
induction, one sees that $T_n$ is a $\pG_n$-stopping time, and the strong
Markov property gives the stated equality in law. Writing
$\Erwt_{x,\omega}$ for the expectation with respect to $\Prwt_{x,\omega}=
\Prw_{x,\omega}\otimes\pQ$, we obtain
\begin{align*}
  \Erw_{x,\omega}\left[\tau_L\right]&= \sum_{z\in
    V_L}\Erwt_{x,\omega}\left[\sum_{n=0}^\infty 1_{\{z\}}(X_n)1_{\{n <
      \tau_L\}}\right]\\
  & =\sum_{z\in
    V_L}\Erwt_{x,\omega}\left[\sum_{n=0}^\infty\sum_{k=T_n}^{T_{n+1}-1}
    1_{\{z\}}(X_k)\right]\\
  &=\sum_{z\in V_L}\Erwt_{x,\omega}\left[\sum_{n=0}^\infty\Big(\sum_{y\in
      V_L}1_{\{y\}}(X_{T_n})\Big)\sum_{k=T_n}^{T_{n+1}-1}
    1_{\{z\}}(X_k)\right]\\
  &=\sum_{y\in
    V_L}\sum_{n=0}^{\infty}\Erwt_{x,\omega}\left[1_{\{y\}}\left(X_{T_n}\right)
    \Erwt_{x,\omega}\left[\sum_{z\in
        V_L}\sum_{k=T_n}^{T_{n+1}-1}1_{\{z\}}(X_k)\bigm|X_{T_n}\right]\right]\\
  &= \sum_{y\in
    V_L}\sum_{n=0}^\infty\Erwt_{x,\omega}\left[1_{\{y\}}\left(X_{T_n}\right)\Erwt_{X_{T_n,\omega}}\left[\sum_{z\in
        V_L}\sum_{k=0}^{T_1-1}1_{\{z\}}\left(X_k\right)\right]\right]\\
  &= \sum_{y\in
    V_L}\sum_{n=0}^\infty\Erwt_{x,\omega}\left[1_{\{y\}}\left(X_{T_n}\right)\right]\Lambda_L(y)
  = \sum_{y\in V_L}\Gh_L(x,y)\Lambda_L(y) = \Gh_L\Lambda_L(x).
\end{align*}
\end{proof1}

\subsection{Space-good/bad and time-good/bad points}
\label{SUBSstgoodandbad}
We classify the grid points inside $V_L$ into good and
bad points, with respect to both space and time. We start by defining
space-good and space-bad points. Unlike in~\cite{BB}, we need
simultaneous control over two scales. This suggests the following stronger
notion of ``goodness''.

\subsubsection{Space-good and space-bad points}
Recall the assignment~\eqref{kernelpL}. We say that $x\in V_L$ is {\it
  space-good} (with respect to $L$ and $\delta >0$), if
\begin{itemize}
\item $x\in V_L\backslash \badP_L$, that is $x$ is good in the sense of
  Section~\ref{SUBSgoodandbad}.
\item If $\dL(x) > 2s_L$, then additionally for all $t\in [h_L(x),2h_L(x)]$
  and for all $y\in V_t(x)$,
  \begin{itemize}
  \item For all $t'\in [h_t^x(y), 2h_t^x(y)]$, with
    $\tilde{q}=p_{h_t^x(y)}$,
    $\left\|(\Pi_{V_{t'}(y)}-\pi_{V_{t'}(y)}^{(\tilde{q})})(y,\cdot)\right\|_1
    \leq \delta$.
  \item If $t-|y-x|  > 2r_t$, then additionally (with the same $\tilde{q}$)
   $$
   \left\|(\Ph_t^x -
     \ph_t^{(\tilde{q}),x})\ph_t^{(\tilde{q}),x}(y,\cdot)\right\|_1 \leq
   (\log h_t^x(y))^{-9}.
   $$
 \end{itemize}
\end{itemize}
In other words, for a point $x\in V_L$ with $\dL(x)>2s_L$ to be space-good,
we do not only require that $x$ is good in the sense of
Section~\ref{SUBSgoodandbad}, but also that all points $y\in V_t(x)$ for
every $t\in[h_L(x),2h_L(x)]$ are good. In this way, we obtain control over
the exit distributions from smaller balls in the bulk of $V_L$, whose radii
are on a preceding scale.

A point $x\in V_L$ which is not space-good is called {\it space-bad}.  The
(random) set of all space-bad points inside $V_L$ is denoted by
$\badPs=\badPs(\omega)$. We write $\goods = \{\badPs = \emptyset\}$ for the
set of environments which contain no bad points, and $\bads=
\left\{\badPs\neq\emptyset\right\}$ for its complement. Notice that in the
notation of Section~\ref{SUBSgoodandbad}, $\badP_L \subset \badPs$ and
$\goods\subset \good$.

On the event $\goods$, we have good control over the coarse grained Green's
functions $\Gh_L$ and  $\Gh_t^x$ in terms of the deterministic kernel
$\Gamma$.
\begin{lemma}
\label{times-superlemma3}
There exist $\delta>0$ small and $L_0$ large such that if $\cspace(\delta,
L_0,L_1)$ is satisfied for some $L_1\geq L_0$, then we have for $L_1\leq
L\leq L_1(\log L_1)^2$ on $\goods$,
\begin{enumerate} 
 \item $\Gh_L \preceq C\Gamma_L$.
 \item If $x\in V_L$ with $\dL(x) > 2s_L$, then for all $t\in[h_L(x),
   2h_L(x)]$,
   $$
  \Gh_t^x \preceq C\Gamma_{t}(\cdot-x,\cdot-x).
$$
\end{enumerate}
\end{lemma}
\begin{proof1}
  (i) Since $\goods\subset\good$, we have $\Gh = \Ghg$ on $\goods$, and
  Lemma~\ref{superlemma} applies. For (ii), with $x$ and $t$ as in
  the statement, there are no bad points within
  $V_t(x)$ on $\goods$. Therefore, also the kernel $\Gh_t^x$ coincides with its
  goodified version, and the
  claim follows again from Lemma~\ref{superlemma}.
\end{proof1}
\begin{lemma}
\label{times-lemmabads}
Assume $\cspace(\delta,L_0,L_1)$. Then for $L$ with $L_1\leq L \leq L_1(\log L_1)^2$,
$$
\pP(\bads) \leq \exp\left(-(2/3)(\log L)^{2}\right).
$$
\end{lemma}
\begin{proof1}
  One can argue as in the proof of~\cite[Lemma 2.3]{BB} (or as in
  the proof of the next Lemma~\ref{times-lemmamanybad}), using
  repeatedly the estimate~\eqref{cspace-essential} under
  $\cspace(\delta,L_0,L_1)$. We omit the details.
\end{proof1}

\subsubsection{Time-good and time-bad points}
We also classify points inside $V_L$ according to the mean time the RWRE
spends in surrounding balls. Remember Condition $\ctime(\eta,L_1)$ and the
function $f_\eta$ introduced above. We now fix $0<\eta<1$.

For points in the bulk of $V_L$, we need control over two scales, as
above. We say that a point $x\in V_L$ is {\it time-good} if the following
holds:
\begin{itemize}
\item For all $x\in V_L$, $t\in[h_L(x),2h_L(x)]$,
$$\Erw_{x,\omega}\left[\tau_{V_{t}(x)}\right]\in
  \left[1-f_\eta(s_L),\, 1+f_\eta(s_L)\right]\cdot
  \Erw_{x}\left[\tau_{V_{t}(x)}\right].
$$  
\item If $\dL(x) > 2s_L$, then additionally for all $t\in
  [h_L(x),2h_L(x)]$, $y\in V_t(x)$ and for all
  $t'\in[h_t^x(y),2h_t^x(y)]$,
$$
\Erw_{y,\omega}\left[\tau_{V_{t'}(y)}\right]\in
  \left[1-f_\eta(s_t),\, 1+f_\eta(s_t)\right]\cdot
  \Erw_{y}\left[\tau_{V_{t'}(y)}\right]. 
$$  
\end{itemize}
A point $x\in V_L$ which is not time-good is called {\it time-bad}.  We
denote by $\badPt = \badPt(\omega)$ the set of all time-bad points inside
$V_L$. With
$$ \mathcal{D}_L = \left\{V_{4h_L(x)}(x) : x\in V_L\right\},$$ 
we let  
$\onebadt= \left\{\badPt\subset D\mbox{ for some }
  D\in\mathcal{D}_L\right\}$, $\manybadt={\left(\onebadt\right)}^c$, and
$\goodt = \{\badPt = \emptyset\}\subset\onebadt$.

The next lemma ensures that for propagating Condition $\ctime$, we can
forget about environments with space-bad points or widely spread time-bad
points.
\begin{lemma}
\label{times-lemmamanybad}
Assume $\cspace(\delta,L_0,L_1)$, and $\ctime(\eta,L_1)$. Then for $L_1\leq L \leq L_1(\log L_1)^2$,
$$
\pP\left(\bads\cup\manybadt\right) \leq
    \exp\left(-(1/2)(\log L)^2\right). 
$$
\end{lemma}
\begin{proof1}
  We have $\pP\left(\bads\cup\manybadt\right) \leq \pP\left(\bads\right) +
  \pP\left(\manybadt\right)$. Lemma~\ref{times-lemmabads} bounds the first
  summand. Now if $x\in\badPt$, then either
  $$
  \Erw_{x,\omega}\left[\tau_{V_t(x)}\right]\notin\left[1-f_\eta(t),\, 1+f_\eta(t)\right]\cdot
  \Erw_x\left[\tau_{V_{t}(x)}\right]
$$
for some $t\in[h_L(x),2h_L(x)]$ (recall that $f_\eta$ is increasing), or, if $\dL(x)>2s_L$, 
  $$
  \Erw_{y,\omega}\left[\tau_{V_{t'}(y)}\right]\notin\left[1-f_\eta(t'),\, 1+f_\eta(t')\right]\cdot
  \Erw_y\left[\tau_{V_{t'}(y)}\right]
$$
for some $y\in V_{2h_L(x)}(x)$, $t'\in[h_{h_L(x)}^x(y),2h_{2h_L(x)}^x(y)]$.

For all $x\in V_L$, we have $h_L(x) \geq r_L/20$. Moreover,
if $\dL(x) > 2s_L$, then $h_L(x)= s_L/20$, whence for all 
$t\in [h_L(x),\, 2h_L(x)]$, all $y\in V_t(x)$, it follows that $h_t^x(y) \geq
r_{(s_L/20)}/20$. Therefore, under $\ctime(\eta,L_1)$,
  $$
\pP\left(x\in\badPt\right) \leq s_L^d\exp\left(-(1/2)\left(\log\left(r_L/20\right)\right)^2\right) +
Cs_L^ds_{s_L}^d\exp\left(-(1/2)\left(\log\left(r_{s_L/20}/20\right)\right)^2\right).
$$
We now observe that if $\omega\in\manybadt$, then there are at least two
time-bad points $x,y$ inside $V_L$ with $|x-y|>2h_L(x)+2h_L(y)$. For such
$x,y$, the events $\{x\in\badPt\}$ and $\{y\in\badPt\}$ are independent.
With the last display, we therefore conclude that
 $$
  \pP\left(\manybadt\right) \leq
  CL^{6d}\left[\exp\left(-(1/2)\left(\log\left(r_{s_L/20}/20\right)\right)^2\right)\right]^2
  \leq \exp\left(-(2/3)(\log L)^2\right). 
  $$
\end{proof1}

\subsection{Proof of the main technical statement}
\label{SUBSctimes-proofs}
In this part, we prove Proposition~\ref{main-prop-times}. We will
always assume that $\delta$ and $L$ are such that
Lemma~\ref{times-superlemma3} can be applied. We start with two auxiliary
statements: Lemma~\ref{times-auxillemma1} proves a difference estimate for mean sojourn
times. Here the difference estimates for the coarse grained Green's
functions from Section~\ref{super-difference} play a crucial
role. Lemma~\ref{times-auxillemma2} then provides the key estimate for
proving the main propagation step.

Note that due to Lemma~\ref{LmeanRWREapriori}, we have for $\omega\in
(\mathcal{P}_\e)^{\mathbb{Z}^d}\cap(\mPsymone)^{\mathbb{Z}^d}$, that is for
$\pP$-almost all environments,
\begin{equation}
\label{generalmeanbound}
\Lambda_L(x) \leq Cs_L^2\leq C(\log L)^{-6}L^2\quad\mbox{for all  } x\in V_L.
\end{equation}

\begin{lemma}
\label{times-auxillemma1}
Assume $\czero(\e)$, $\cB$, $\cspace(\delta, L_0,L_1)$, and let $L_1\leq L\leq L_1(\log L_1)^2$.
Let $0\leq \alpha < 3$ and $x, y\in
V_{L-2s_L}$ with $|x-y| \leq (\log
s_L)^{-\alpha}\,s_L$. Then for $\pP$-almost all $\omega\in \goods$, 
$$
\left|\Lambda_L(x) - \Lambda_L(y)\right| \leq C(\log\log s_L)(\log s_L)^{-\alpha}
  s_L^2.
$$
\end{lemma}
\begin{proof1}
We let
$\omega\in(\mathcal{P}_\e)^{\mathbb{Z}^d}\cap(\mPsymone)^{\mathbb{Z}^d}\cap\goods$.
The statement follows if we show that for all $t\in
\left[(1/20)s_L,(1/10)s_L\right]$,
$$
  \left|\Erw_{x,\omega}\left[\tau_{V_t(x)}\right] -
    \Erw_{y,\omega}\left[\tau_{V_t(y)}\right]\right| \leq C(\log\log t)(\log t)^{-\alpha}
  t^2.
$$
We fix such a $t$ and set $t' = \left(1-20(\log t)^{-\alpha}\right)t$. Then $V_{t'}(x) \subset
  V_t(x)\cap V_t(y)$. Now put $B = V_{t'-2s_t}(x)$. By
  Lemma~\ref{times-keylemma}, we have the representation 
\begin{equation}
\label{times-auxillemma1-1}
\Erw_{x,\omega}\left[\tau_{V_t(x)}\right] = \Gh_t^x1_B\Lambda_t^x(x) +
 \Gh_t^x1_{V_{t}(x)\backslash B}\Lambda_t^x(x).
\end{equation}
By~\eqref{generalmeanbound}, $\Lambda_t^x(z)\leq C(\log t)^{-6}t^2$, for
all $z\in V_t(x)$. Moreover, since $\omega \in \goods$, we have by
Lemma~\ref{times-superlemma3} $\Gh_t^x\preceq
C\Gamma_t(\cdot-x,\cdot-x)$. Applying Lemma~\ref{super-gammalemma} (ii)
gives
$$
\Gh_t^x1_{V_{t}(x)\backslash B}\Lambda_t^x(x)\leq
C\Gamma_t\left(0,V_t\backslash V_{t'-2s_t}\right)\,(\log t)^{-6}t^2 \leq
(\log t)^{-\alpha}t^2
$$ 
for $L$ (and therefore also $t$) sufficiently large. 
Concerning $\Erw_{y,\omega}\left[\tau_{V_t(y)}\right]$, we write again
$$
  \Erw_{y,\omega}\left[\tau_{V_t(y)}\right]= \Gh_t^y1_B\Lambda_t^y(y) + \Gh_t^y1_{V_t(y)\backslash B}\Lambda_t^y(y).
$$
The second summand is bounded by $(\log t)^{-\alpha}t^2$, as in the display
above. For $z \in B$, we have $h_t^x(z) = h_t^y(z) = (1/20)s_t$. In
particular, $\Ph_t^x(z,\cdot) = \Ph_t^y(z,\cdot)$, and also $\Lambda_t^x(z)
= \Lambda_t^y(z)$. Since both $x$ and $y$ are contained in $B\subset
V_t(x)\cap V_t(y)$, the strong Markov property gives
$$
\Gh_t^y(y,z) =
\Gh_t^x(y,z) +b(y,z),
$$
where
$$
  b(y,z) = \Erw_{y,\Ph_t^y (\omega)}\left[\Gh_t^y(\tau_B,z);\,
    \tau_B < \infty\right] -\Erw_{y,\Ph_t^x(\omega)}\left[\Gh_t^x(\tau_B,z);\,
    \tau_B < \infty\right].
$$
Therefore,
\begin{align*}
  \lefteqn{\left|\Erw_{x,\omega}\left[\tau_{V_t(x)}\right]-\Erw_{y,\omega}\left[\tau_{V_t(y)}\right]\right|}\\
  &\leq 2(\log t)^{-\alpha}t^2 +\sum_{z\in B}\left(\big|\Gh_t^x(x,z)-\Gh_t^x(y,z)\big| + |b(y,z)|\right)\Lambda_t^x(z).
\end{align*}
The quantity $\Lambda_t^x(z)$ is again bounded by $C(\log
t)^{-6}t^2$. For the part of the sum involving $|b(y,z)|$, we
notice that if $w\in V_t(y)\backslash B$, then $t-|w-y|\leq C(\log
t)^{-\alpha}t$ and similarly for $v\in V_t(x)$. We can
use twice Lemma~\ref{super-gammalemma} (ii) to get
$$
\sum_{z\in B}|b(y,z)| \leq \sup_{v\in
  V_t(x)\backslash B}\Gh_t^x(v,B) + \sup_{w\in
  V_t(y)\backslash B}\Gh_t^y(w,B) \leq C(\log t)^{6-\alpha}.
$$
Finally, for the sum over the Green's function difference, we recall that
$\Gh_{t}^x$ coincides with its goodified version. Applying 
Lemma~\ref{super-greendifference} $O\left((\log
  t)^{3-\alpha}\right)$ times gives
$$
\sum_{z\in B}\big|\Gh_{t}^x(x,z)-\Gh_{t}^x(y,z)\big| \leq
C(\log\log t)(\log t)^{6-\alpha}.
$$
This proves the statement.
\end{proof1}

\begin{lemma}
\label{times-auxillemma2}
Assume $\czero(\e)$, $\cB$, $\cspace(\delta, L_0,L_1)$, and let $L_1\leq L\leq
L_1(\log L_1)^2$. Let $p=p_{s_L/20}$, cf.~\eqref{kernelpL}, and set $\Delta
= 1_{V_L}(\Ph_L-\ph_L^{(p)})$. For $\pP$-almost all $\omega\in\goods$,
$$
\sup_{x\in V_L}\big|\Gh_L\Delta\gh_L^{(p)}\Lambda_L(x)\big| \leq C(\log L)^{-5/3}L^2.
$$ 
\end{lemma}
\begin{proof1}
  Again, we consider
  $\omega\in(\mathcal{P}_\e)^{\mathbb{Z}^d}\cap(\mPsymone)^{\mathbb{Z}^d}\cap\goods$.
  Write $\gh= \gh^{(p)}$, $\ph=\ph^{(p)}$. First,
$$
\Gh\Delta\gh\Lambda_L(x) = \Gh\Delta\ph\gh\Lambda_L(x) + 
\Gh\Delta\Lambda_L(x)= A_1 + A_2.
$$
By Lemma~\ref{times-superlemma3}, $\Gh = \Ghg \preceq
C\Gamma$. Therefore, with $B_1= V_{L-2r_L}$, we bound $A_1$ by
\begin{align*}
  \left|A_1\right| &\leq \big|\Gh1_{B_1}\Delta\ph\gh\Lambda_L(x)\big| +
  \big|\Gh1_{V_L\backslash B_1}\Delta\ph\gh\Lambda_L(x)\big|\\
  &\leq \Big|\sum_{v\in B_1, w\in V_L}\Gh(x,v)\Delta\ph(v,w)\sum_{y\in
    V_L}\left(\gh(w,y)-\gh(v,y)\right)\Lambda_L(y)\Big|+ C(\log L)^{-2}L^2 \\
  &\leq C(\log L)^{-5/3}L^2,
 \end{align*}
 where in the next to last inequality we have used the bound on
 $\Lambda_L(y)$ given by~\eqref{generalmeanbound} and
 Lemma~\ref{super-gammalemma} (ii), (iii) for $\Gamma$, and in the last inequality
 additionally Lemma~\ref{super-greendifference}.  For the term $A_2$, we
 let $B= V_{L-5s_L}$ and split into
$$
  A_2 = \Gh1_{B}\Delta\Lambda_L(x) + \Gh1_{V_L\backslash B}\Delta\Lambda_L(x).
$$
Lemma~\ref{super-gammalemma} (ii) yields
$$
\Gh(x,\sh_L(5s_L))\leq
C\log\log L.
$$
Since $\Lambda_L(y)\leq (\log L)^{-2}L^2$ by~\eqref{generalmeanbound}, this
is good enough for the second summand of $A_2$. For the first one,
$$
\Gh1_B\Delta\Lambda_L(x) \leq C\Gamma(x,B)\sup_{v\in
  B}\left|\Delta\Lambda_L(v)\right|.
$$
Since $\Gamma(x,B) \leq C(\log L)^6$, the claim follows we show that for
$v\in B$,  
\begin{equation}
\label{times-mainlemma-1}
\left|\Delta\Lambda_L(v)\right| \leq C(\log L)^{-8}L^2,
\end{equation}
which, by definition of $\Delta$, in turn follows if for all
$t\in[h_L(v), 2h_L(v)]$,
$$
\left|\left(\Pi_{V_t(v)}-\pi_{V_t(v)}\right)\Lambda_L(v)\right| \leq C(\log L)^{-8}L^2,
$$
where we have set $\pi_{V_t(v)} = \pi_{V_t(v)}^{(p)}$. Notice that on $B$,
$h_L(\cdot) = (1/20)s_L$. We now fix $v\in B$ and $t\in[(1/20)s_L,
(1/10)s_L]$. Set $\Delta'= 1_{V_t(v)}(\Ph_t^v-\ph_t^{(p),v})$ and $B'=
V_{t-2r_t}(v)$. By expansion~\eqref{prel-pbe1},
\begin{equation}
\label{times-mainlemma-2}
\left(\Pi_{V_t(v)}-\pi_{V_t(v)}\right)\Lambda_L(v) =
\Gh_t^v1_{B'}\Delta'\pi_{V_t(v)}\Lambda_L(v) +
\Gh_t^v1_{V_t(v)\backslash B'}\Delta'\pi_{V_t(v)}\Lambda_L(v).
\end{equation}
Since $\pi_{V_t(v)} = \ph_t^{(p),v}\pi_{V_t(v)}$, we get
\begin{align*}
\big|\Gh_t^v1_{B'}\Delta'\pi_{V_t(v)}\Lambda_L(v)\big| &\leq 
  \Gh_t^v(v,B')\sup_{w\in
    B'}\big|\Delta'\ph_t^v(w,\cdot)\big|_1\sup_{y\in \partial V_t(v)}\Lambda_L(y)\\
&\leq C(\log s_L)^6(\log L)^{-6}L^2\sup_{w\in
    B'}\big\|\Delta'\ph_t^{(p),v}(w,\cdot)\big\|_1.
\end{align*}
Here, in the last inequality we have used~\eqref{generalmeanbound} and
Lemma~\ref{super-gammalemma} (iii). In order to bound
$\|\Delta'\ph_t^{(p),v}(w,\cdot)\|_1$ for $w\in B'$, we use the fact that
$v$ is space-good and $\dL(v)>2s_L$, which gives also control over the exit
distributions from smaller balls inside $V_t(v)$. Indeed, by definition we first have for
$w\in B'$, with $\tilde{q}=p_{h_t^x(y)}$,
$$\big\|1_{V_t(v)}(\Ph_t^v-\ph_t^{(\tilde{q}),v})\ph_t^{(\tilde{q}),v}(w,\cdot)\big\|_1
\leq (\log h_t^v(w))^{-9}\leq C(\log L)^{-9}.$$ The last inequality follows
from the bound $h_t^v(w) \geq (1/20)r_{s_L/20}$.  Furthermore, under
$\cspace(\delta,L_0,L_1)$, the kernel $\tilde{q}$ is close to $p$: in fact,
one has $\|\tilde{q}-p\|_1\leq C(\log L)^{-8}$, see~\cite[Lemma 2.2]{BB}
and the arguments in the proof there. This bound transfers to the exit
measures, so that, arguing as in~\cite[Lemma 4.1]{BB},
$$
 \sup_{w\in B'}\big\|\Delta'\ph_t^{(p),v}(w,\cdot)\big\|_1 =\sup_{w\in
   B'}\big\|(\Ph_t^v-\ph_t^{(p),v})\ph_t^{(p),v}(w,\cdot)\big\|_1 \leq C(\log L)^{-8}.
$$
Putting the estimates together, we obtain as desired
$$
\big|\Gh_t^v1_{B'}\Delta'\pi_{V_t(v)}\Lambda_L(v)\big| \leq C(\log L)^{-8}L^2.
$$
For the second summand of~\eqref{times-mainlemma-2},
Lemma~\ref{super-gammalemma} (ii) gives $\Gh_t^v(v,V_t(v)\backslash B')
\leq C$, whence
$$
\big|\Gh_t^v1_{V_t(v)\backslash
  B'}\Delta'\pi_{V_t(v)}\Lambda_L(v)\big| \leq C \sup_{w\in V_t(v)\backslash
  B'}\big|\Delta'\pi_{V_t(v)}\Lambda_L(w)\big|. 
$$
Fix $w\in V_t(v)\backslash B'$. Set $\eta = \dist(w,\partial V_t(v))\leq
2r_t + \sqrt{d}$ and choose $y_w\in\partial V_t(v)$ such that $|w-y_w| =
\eta$. With
$$
I(y_w) = \left\{y\in\partial
    V_t(v) : |y-y_w|\leq (\log L)^{-5/2}s_L\right\},
$$
we write
\begin{align}
\label{times-mainlemma-3}
\Delta'\pi_{V_t(v)}\Lambda_L(w)& = \sum_{y\in\partial
  V_t(v)}\Delta'\pi_{V_t(v)}(w,y)\left(\Lambda_L(y)-\Lambda_L(y_w)\right)\nonumber\\
&= \sum_{y \in
  I(y_w)}\Delta'\pi_{V_t(v)}(w,y)\left(\Lambda_L(y)-\Lambda_L(y_w)\right)\nonumber\\
&\quad +
\sum_{y\in\partial
  V_t(v)\backslash
  I(y_w)}\Delta'\pi_{V_t(v)}(w,y)\left(\Lambda_L(y)-\Lambda_L(y_w)\right).
 \end{align}
For $y\in I(y_w)$, Lemma~\ref{times-auxillemma1} yields
$|\Lambda_L(y)-\Lambda_L(y_w)| \leq C(\log L)^{-7/3}s_L^2$. Therefore, 
$$
\sum_{y \in
  I(y_w)}\left|\Delta'\pi_{V_t(v)}(w,y)\right|\left|\Lambda_L(y)-\Lambda_L(y_w)\right|
\leq C(\log L)^{-8}L^2.
$$
It remains to handle the second term of~\eqref{times-mainlemma-3}. To
this end, let $U(w) = \{u\in V_t(v) : |\Delta'(w,u)|> 0\}$.
Using for $y\in \partial V_t(v)\backslash I(y_w)$ the trivial bound
$$\left|\Lambda_L(y)-\Lambda_L(y_w)\right|\leq \Lambda_L(y) + \Lambda_L(y_w)
\leq C(\log L)^{-6}L^2,$$
see~\eqref{generalmeanbound}, we obtain
\begin{multline*}
\sum_{y\in\partial
    V_t(v)\backslash
    I(y_w)}\left|\Delta'\pi_{V_t(v)}(w,y)\left(\Lambda_L(y)-\Lambda_L(y_w)\right)\right|\\
\leq C(\log L)^{-6}L^2\sup_{u\in U(w)}\pi_{V_t(v)}\left(u,\partial
V_t(v)\backslash I(y_w)\right).
\end{multline*}
If $u\in U(w)$ and $y\in\partial V_t(v)\backslash I(y_w)$, then 
$$
|u-y| \geq |y-y_w| -|y_w-u| \geq (\log L)^{-5/2}s_L - 3r_t \geq (1/2)(\log L)^{-5/2}s_L.
$$
For such $u$, we get by standard hitting estimates, see e.g.~\cite[Lemma
3.2 (ii)]{BB}, 
\begin{align*}
\pi_{V_t(v)}\left(u,\partial V_t(v)\backslash I(y_w)\right) &\leq C r_t
\sum_{y\in \partial V_t(v)\backslash I(y_w)}\frac{1}{|u-y|^d} \\
& \leq Cr_t(\log L)^{5/2}(s_L)^{-1} \leq C(\log L)^{-9}.
\end{align*}
The estimate on the sum can be obtained from~\cite[Lemma
3.6]{BB}. This bounds the second term of~\eqref{times-mainlemma-3}. We have
proven~\eqref{times-mainlemma-1} and hence the lemma.
\end{proof1}

Now it is easy to prove the main propagation step.
\begin{lemma}
\label{times-mainlemma}
Assume $\czero(\e)$ and $\cB$. 
There exists $L_0=L_0(\eta)$ such that if $L_1\geq L_0$ and 
$\cspace(\delta,L_0,L_1)$ holds, then for $L_1\leq L\leq L_1(\log L_1)^2$ and
$\pP$-almost all $\omega\in \goods\cap\onebadt$,
$$
  \Erw_{0,\omega}\left[\tau_L\right] \in
    \left[1-f_\eta(L),\,1+f_\eta(L)\right]\cdot\Erw_0\left[\tau_L\right].
$$
\end{lemma}
\begin{proof1}
We let $\omega\in(\mathcal{P}_\e)^{\mathbb{Z}^d}\cap(\mPsymone)^{\mathbb{Z}^d}\cap\goods\cap\onebadt$.
 Put $p=p_{s_L/20}$. In this proof, we keep the superscript $(p)$ in $\ghp$.
  By Lemma~\ref{times-keylemma} and the perturbation
  expansion~\eqref{prel-pbe1}, with $\Delta = 1_{V_L}(\Ph-\php)$,
$$
  \Erw_{0,\omega}\left[\tau_L\right] = \Gh\Lambda_L(0)
  = \ghp\Lambda_L(0) +\Gh\Delta\ghp\Lambda_L(0) = A_1 + A_2.
$$
Set $B= V_L\backslash (\badPt\cup\sh_L(L/(\log L)^2)$. The term $A_1$ we split into
$$
A_1 = \ghp1_B\Lambda_L(0) + \ghp1_{V_L\backslash B}\Lambda_L(0).
$$
Since $\ghp(0,V_L\backslash B) \leq C(\log L)^3$ by Lemma~\ref{super-gammalemma} (ii)
and $\Lambda_L(x) \leq (\log
L)^{-6}L^2$, the second summand of $A_1$ can be bounded by $O\left((\log
  L)^{-3}\right)\Erw_0[\tau_L]$. The main contribution comes from the first
summand. First notice that
$$
\ghp1_B\lambdap_L(0) =
\ghp1_B\lambda_L(0)\left(1+O\left(s_L^{-1}\right)\right)
=\Erw_{0}\left[\tau_L\right]\left(1+O\left((\log L)^{-6}\right)\right).
$$
Furthermore, we have for $x\in B$ by definition
$$
  \Lambda_L(x) \in \left[1-f_\eta\left((\log L)^{-3}L\right),\, 1+f_\eta\left((\log L)^{-3}L\right)\right]\cdot\lambda_L(x).
$$
Collecting all terms, we conclude that
\begin{equation*}
\begin{split}
  A_1 \in &\left[1 - O\left((\log L)^{-3}\right)- f_\eta\left((\log
      L)^{-3}L\right),\, 1 + O\left((\log L)^{-3}\right) + f_\eta\left((\log
      L)^{-3}L\right)\right]\\
&\times\Erw_0\left[\tau_L\right]. 
\end{split}
\end{equation*}
Lemma~\ref{times-auxillemma2} bounds $A_2$ by $O((\log
L)^{-5/3})\Erw_0[\tau_L]$. Since for $L$ sufficiently large,
$$
f_\eta(L) > f_\eta\left((\log L)^{-3}L\right) +C(\log L)^{-5/3},
$$
we arrive at
$$
  \Erw_{0,\omega}\left[\tau_L\right] = A_1 + A_2 \,\in\, 
  \left[1-f_\eta(L),\,1+f_\eta(L)\right] \cdot\Erw_0\left[\tau_L\right],
$$
as claimed.
\end{proof1}
Proposition~\ref{main-prop-times} follows now as an immediate
consequence of our estimates.

\begin{proof2}{\bf of Proposition~\ref{main-prop-times}:}
  (i) From Lemmata~\ref{times-lemmamanybad}
  and~\ref{times-mainlemma} we deduce that for large $L_0$, if $L_1\geq
  L_0$ and $L_1\leq L\leq L_1(\log L_1)^2$, we have under
  $\cspace(\delta,L_0,L_1)$ and $\ctime(\eta,L_1)$
  \begin{align*}
    \lefteqn{\pP\left(\Erw_{0,\omega}\left[\tau_L\right] \notin
        \left[1-f(L),\,
          1+f(L)\right]\cdot\Erw_0\left[\tau_L\right]\right)}\\
    &\leq \pP\left(\bads\cup\manybadt\right)+ \\
    &\quad +\; \pP\left(\Erw_{0,\omega}\left[\tau_L\right] \notin
        \left[1-f(L),\,1+f(L)\right]\cdot\Erw_0[\tau_L];\,\goods\cap\onebadt\right)\\
    &\leq\exp\left(-(1/2)(\log L)^2\right).
  \end{align*}
  By Proposition~\ref{main-prop-exitmeas}, if $\delta>0$ is small and $L_0$ is
  sufficiently large, $\cspace(\delta,L_0,L)$ holds under $\czero(\e)$ for all
  $L\geq L_0$, provided
  $\e\leq\e_0(\delta)$. This proves the proposition.
\end{proof2}

\subsection{Proof of the main theorem on sojourn times}
We shall first prove convergence of the (non-random) sequence
$\pE[\Erw_{0,\omega}[\tau_L]]/L^2$ towards a constant $D$ that lies in a
small interval around $1$. Note that Proposition~\ref{main-prop-times}
together with Lemma~\ref{LmeanRWREapriori} already tells us that for any
$0<\eta<1$, under $\czero$, $\cB$ and $\cone(\e)$ for $\e(\eta)$ small,
$$
\pE\left[\Erw_{0,\omega}\left[\tau_L\right]\right]/L^2
\in[1-\eta,1+\eta]\quad\hbox{for large }L.
$$

\begin{proposition}
\label{prop-times-underP}
Assume $\cone$ and $\cB$. Given $0<\eta < 1$, one can find $\e_0=\e_0(\eta)
> 0$ such that if $\czero(\e)$ is satisfied for some $\e\leq \e_0$, then
there exists $D\in [1-\eta, 1+\eta]$ such that
$$
\lim_{L\rightarrow\infty}\left(\pE\left[\Erw_{0,\omega}\left[\tau_L\right]\right]/L^2\right)=D.
$$ 
\end{proposition}
\begin{proof1}
  Let $0<\eta<1$. By choosing first $\delta$, then $L_0$ and then $\e_0$
  small respectively large enough, we know from
  Propositions~\ref{main-prop-exitmeas} and~\ref{main-prop-times} that
  under $\cone$ and $\cB$, whenever $\czero(\e)$ is satisfied for some
  $\e\leq \e_0$, $\cspace(\delta,L_0,L)$ and $\ctime(\eta/2,L)$ hold true
  for all $L\geq L_0$. We can therefore assume both conditions. We obtain
  from Lemma~\ref{times-lemmabads}
$$
\pE\left[\Erw_{0,\omega}\left[\tau_L\right]\right] =
\pE\left[\Erw_{0,\omega}\left[\tau_L\right];\,\goods\right]
+O\left(L^2\exp\left(-(2/3)(\log L)^2\right)\right).
$$
Thus it suffices to look at $\Erw_{0,\omega}[\tau_L]$ on the event
$(\mathcal{P}_\e)^{\mathbb{Z}^d}\cap(\mPsymone)^{\mathbb{Z}^d}\cap\goods$. Setting $B= V_L\backslash (\sh_L(L/(\log
L)^2))$, we see from the proof of Lemma~\ref{times-mainlemma} 
that on this this event,
$$
\Erw_{0,\omega}\left[\tau_L\right] = (\ghp1_B\Lambda_L)(0) + O\left((\log L)^{-{5/3}}L^2\right),
$$
where the constant in the error term does only depend on $d$ (and not on
$L$ or the environment). For $x\in B$, $h_L(x)= s_L/20$. In particular, this implies on the set $B$
$$
\pE\left[\Lambda_L(\cdot)\right] \equiv
\pE\left[\Lambda_L(0)\right],\quad\hbox{and } \lambdap(\cdot) \equiv \lambdap(0).
$$
Now put $c_L = \pE\left[\Lambda_L(0)\right]/\lambdap(0)$. We
have
\begin{align*}
\pE\left[\Erw_{0,\omega}\left[\tau_L\right]\right]
&= \ghp(0,B)\cdot\pE\left[\Lambda_L(0)\right] + O\left((\log L)^{-5/3}L^2\right)\\
&= c_L\,\ghp(0,B)\cdot\lambdap(0) +  O\left((\log
  L)^{-5/3}L^2\right)\\
&= c_L\,\Erw_{0,p}\left[\tau_L\right] +O\left((\log
  L)^{-5/3}L^2\right).
\end{align*}
Since $\Erw_{0,p}\left[\tau_L\right]/L^2$ converges to $1$ by
Lemma~\ref{LmeanRWp}, convergence of
$\pE\left[\Erw_{0,\omega}[\tau_L]\right]/L^2$ follows if we show that 
$\lim_{L\rightarrow\infty}c_L$ exists. Let $L'\in (L,2L]$. As before,
\begin{equation}
\label{eq-obvious}
\pE\left[\Erw_{0,\omega}[\tau_{L'}]\right]=c_{L'}\,\Erw_{0,p}\left[\tau_{L'}\right] +O\left((\log
  L)^{-5/3}L^2\right).
\end{equation}
On the other hand, we claim that~\eqref{eq-obvious} also holds with $c_{L'}$
replaced by $c_L$. To see this, we slightly change the coarse graining
scheme inside $V_{L'}$, as in the proof of~\cite[Proposition 1.1]{BB}.
More specifically, we define for $L'\in (L,2L]$ the
coarse graining function 
$\tilde{h}_{L'} : \overline{C}_{L'} \rightarrow\mathbb{R}_+$ by setting 
$$
 \tilde{h}_{L'}(x) = \frac{1}{20}\max\left\{s_L
  h\left(\frac{\textup{d}_{L'}(x)}{s_{L'}}\right),\, r_L\right\}.
$$
Then $\tilde{h}_{L'}(x) = h_L(0)=s_L/20$ for $x\in V_{L'}$ with
$\textup{d}_{L'}(x) \geq 2s_{L'}$. We consider the analogous definition of
space-good/bad and time-good/bad points within $V_{L'}$, which uses the
coarse graining function $\tilde{h}_{L'}$ instead of $h_{L'}$ and the
coarse grained transition kernels $\tilde{\Pi}$ and $\tilde{\pi}$ in
$V_{L'}$ defined in terms of $\tilde{h}_{L',r}$,
cf.~\eqref{eq-cgkernel}. Clearly, all the above statements of this
section remain true (at most the constants change), and we can work
with the same kernel $p=p_{s_L/20}$.  Denoting by $\gtp$ the Green's
function corresponding to $\ptp$ and by $\tilde{B}$ the set
$V_{L'}\backslash \sh_{L'}(L'/(\log L')^2)$, we obtain as above
\begin{align*}
\pE\left[\Erw_{0,\omega}[\tau_{L'}]\right]
&= \gtp(0,\tilde{B})\pE\left[\Lambda_L(0)\right] + O\left((\log L)^{-5/3}L^2\right)\\
&= c_L\,\gtp(0,B')\lambdap(0) +  O\left((\log
  L)^{-5/3}L^2\right)\\
&= c_L\,\Erw_{0,p}\left[\tau_{L'}\right] +O\left((\log
  L)^{-5/3}L^2\right).
\end{align*}
Note that since $\tilde{h}_{L'}(\cdot)\equiv h_L(0)$ on $B'$, the
quantities $c_L$, $\pE\left[\Lambda_L(0)\right]$ and
$\lambdap(0)$ do indeed appear in the above display.  Comparing with~\eqref{eq-obvious}, this shows
that for some constant $C>0$
$$|c_L-c_{L'}|\leq C(\log L)^{-5/3},$$
which readily implies that $c_L$ is a Cauchy sequence and thus
$\lim_{L\rightarrow\infty}c_L=D$ exists.  From
Proposition~\ref{main-prop-times} we already know that
$D\in[1-\eta,1+\eta]$. This finishes the proof.
\end{proof1}

We shall now employ Hoeffding's inequality to show that
$\Erw_{0,\omega}[\tau_L]$ is close to its mean.
\begin{lemma}
\label{Ltimes-concentration}
Assume $\cone$ and $\cB$. There exists $\e_0 >
0$ such that if $\czero(\e)$ holds for some $\e\leq \e_0$, then
$$
\pP\left(\frac{1}{L^2}\big|\Erw_{0,\omega}\left[\tau_L\right]-\pE\left[\Erw_{0,\omega}\left[\tau_L\right]\right]\big|>(\log
  L)^{-4/3}\right)\leq \exp\left(-(1/3)(\log L)^2\right).
$$
\end{lemma}
Let us first show how to prove Theorem~\ref{thm-times} from this result.

\begin{proof2}{\bf of Theorem~\ref{thm-times}:}
We know from Proposition~\ref{prop-times-underP}, $D$ the
constant from there, 
$$
\Big|\frac{1}{L^2}\Erw_{0,\omega}\left[\tau_L\right]-D\Big| \leq
\frac{1}{L^2}\big|\Erw_{0,\omega}\left[\tau_L\right]-\pE\left[\Erw_{0,\omega}\left[\tau_L\right]\right]\big|
+ \alpha(L)
$$
for some (deterministic) sequence $\alpha(L)\rightarrow 0$ as
$L\rightarrow\infty$. Putting 
$$\alpha'(L) = \max\left\{ (\log L)^{-4/3}, \alpha(L)\right\},$$ we deduce from
Lemma~\ref{Ltimes-concentration} that
$$
\pP\left(\Big|\frac{1}{L^2}\Erw_{0,\omega}\left[\tau_L\right]-D\Big|\geq
  2\alpha'(L)\right)\leq \exp\left(-(1/3)(\log L)^2\right).
$$
This implies the first statement of the theorem. For the second, we have
\begin{align*}
\pP\left(\Big|\sup_{x: |x| \leq
      L^k}\Erw_{x,\omega}\left[\tau_{V_L(x)}\right]/L^2 -D\Big| \geq 2\alpha'(L)\right)
&\leq CL^{kd} \pP\left(\left|\Erw_{0,\omega}\left[\tau_L\right]/L^2
    -D\right| \geq 2\alpha'(L)\right)\\
&\leq \exp\left(-(1/4)(\log L)^2\right)
\end{align*}
for large $L$, and the same bound holds with the supremum over $x$ with
$|x|\leq L^k$ replaced by the infimum.  The second claim of the theorem
follows now from Borel-Cantelli.
\end{proof2}

It remains to prove Lemma~\ref{Ltimes-concentration}.

\begin{proof2}{\bf of Lemma~\ref{Ltimes-concentration}:}
By Proposition~\ref{main-prop-exitmeas} and 
Lemma~\ref{times-lemmabads}, we can find $\e_0>0$ such that under
$\czero$ and $\cone(\e)$ for $\e\leq \e_0$,
$$
\pP\left(\bads\right) \leq
    \exp\left(-(2/3)(\log L)^2\right)\quad\hbox{for }L\hbox{ large.}
$$
As in the proof of
Proposition~\ref{prop-times-underP} (or Lemma~\ref{times-mainlemma}), we
have for $\omega\in
(\mathcal{P}_\e)^{\mathbb{Z}^d}\cap(\mPsymone)^{\mathbb{Z}^d}$ in the 
complement of $\bads$, that is $\pP$-almost surely on the event $\goods$,
$$
\Erw_{0,\omega}\left[\tau_L\right] = (\ghp 1_B\Lambda_L)(0) + O\left((\log
  L)^{-{5/3}}L^2\right),
$$
where $B= V_L\backslash (\sh_L(L/(\log
L)^2))$. In the proof of
Proposition~\ref{prop-times-underP} we have also seen that
$$
\pE\left[\Erw_{0,\omega}\left[\tau_L\right]\right] = \ghp(0,B)\pE\left[\Lambda_L(0)\right] + O\left((\log L)^{-5/3}L^2\right).
$$
Therefore, on $\goods$,
\begin{align*}
  \Erw_{0,\omega}\left[\tau_L\right]
  &=
  \ghp(0,B)\pE\left[\Lambda_L(0)\right] +
  \sum_{y\in B}\ghp(0,y)\left(\Lambda_L(y)-\pE\left[\Lambda_L(0)\right]\right) + O\left((\log
    L)^{-5/3}L^2\right)\\
  &= \pE\left[\Erw_{0,\omega}\left[\tau_L\right]\right] + \sum_{y\in B}\ghp(0,y)\left(\Lambda_L(y)-\pE\left[\Lambda_L(0)\right]\right) + O\left((\log
    L)^{-5/3}L^2\right).
\end{align*}
The statement of the lemma will thus follow if we show that
\begin{equation}
\label{eq1-concentration}
\pP\left(\Big|\sum_{y\in B}\ghp(0,y)\left(\Lambda_L(y)-\pE\left[\Lambda_L(0)\right]\right)\Big|>(\log
  L)^{-3/2}L^2\right)\leq \exp\left(-(\log L)^2\right).
\end{equation}
We use a similar strategy as in the proof of~\cite[Lemma 5.4]{BB}. First,
define for $j\in\mathbb{Z}$ the interval $I_j = (js_L,(j+1)s_L]$. Now divide $B$ into subsets $W_{\bj} = B\cap
\left(I_{j_1}\times\dots\times I_{j_d}\right)$, where ${\bj} =
(j_1,\dots,j_d) \in\mathbb{Z}^d$.  Let $J$ be the set of those $\bj$ for
which $W_{\bj} \neq \emptyset$. Then there exists a constant $K=K(d)$ and a
disjoint partition of $J$ into sets $J_1,\dots,J_K$, such that for any
$1\leq \ell\leq K$,
\begin{equation}
\label{eq2-concentration}
\bj,\bj' \in J_\ell,\ \bj \neq \bj' \Longrightarrow
\dist(W_{\bj},W_{\bj'}) > s_L.
\end{equation}
We set 
$$
\xi_{\bj} =  \sum_{y\in W_{\bj}}\ghp(0,y)\left(\Lambda_L(y)-\pE\left[\Lambda_L(0)\right]\right)
$$
and $t = t(d,L) = (\log L)^{-3/2}L^2$.
From~\eqref{eq2-concentration} we see that the random variables
$\xi_{\bj}$, $\bj\in J_\ell$, are independent and centered (we recall again that
$\pE\left[\Lambda_L(y)\right]= \pE\left[\Lambda_L(0)\right]$ for $y\in
B$). Put $\Omega'=(\mathcal{P}_\e)^{\mathbb{Z}^d}\cap(\mPsymone)^{\mathbb{Z}^d}$.
Applying Hoeffding's inequality, we obtain with
${\|\xi_{\bj}\|}_{\infty} = \sup_{\omega\in\Omega'}|\xi_{\bj}(\omega)|$,
for some constant $c>0$,
\begin{equation}
\label{eq3-concentration}
  \pP\left(\Big|\sum_{\bj\in J}\xi_{\bj}\Big| > t\right)\leq 
 K \max_{1\leq \ell\leq K}\pP\left(\Big|\sum_{\bj\in
        J_\ell}\xi_{\bj}\Big| > \frac{t}{K}\right)\leq 2\exp\left(-c\frac{(\log L)^{-3}L^{4}}{\sum_{\bj \in
      J_\ell}{\|\xi_{\bj}\|}^2_{\infty}}\right). 
\end{equation}
It remains to estimate the $\sup$-norm of the
$\xi_{\bj}$. We have, by Lemmata~\ref{superlemma} and~\ref{super-gammalemma},
$$ \gh(x,W_{\bj}) \leq 
    \frac{Cs_L^d}{s_L^2(s_L + \dist(x,W_{\bj}))^{d-2}} =
    C\left(1+\frac{\dist(x,W_{\bj})}{s_L}\right)^{2-d}.
$$
By~\eqref{generalmeanbound},
$$
\left|\Lambda_L(y)-\pE\left[\Lambda_L(0)\right]\right|\leq C
(\log L)^{-6}L^2.
$$
Altogether, recalling that $d\geq 3$,
$$
\sum_{\bj \in
      J_\ell}{\|\xi_{\bj}\|}^2_{\infty} \leq C\sum_{r=1}^{C(\log
      L)^3}r^{-d+3}L^4(\log L)^{-12}\leq C(\log L)^{-9}L^4.
$$
Going back to~\eqref{eq3-concentration}, this shows
$$
\pP\left(\Big|\sum_{y\in B}\ghp(0,y)\left(\Lambda_L(y)-\pE\left[\Lambda_L(0)\right]\right)\Big|\geq(\log
  L)^{-3/2}L^2\right)\leq 2\exp\left(-c(\log L)^{6}\right), 
$$
which is more than we need, cf.~\eqref{eq1-concentration}. This completes the proof of the lemma.
\end{proof2}

\begin{proof2}{\bf of Corollary~\ref{cor-times-moments}:}
  Let $k\in\mathbb{N}$, and let first $m=1$. By 
  Proposition~\ref{main-prop-times}, we obtain under our conditions (for
  $\e$ small) 
\begin{multline*}
\pP\left(\sup_{x: |x| \leq L^{k}}\sup_{y\in
    V_L(x)}\Erw_{y,\omega}\left[\tau_{V_L(x)}\right]/L^2\geq 2\right)
\leq CL^{kd}\pP\left(\sup_{y\in
    V_L}\Erw_{y,\omega}\left[\tau_L\right]/L^2\geq 2\right)\\
\leq CL^{(k+1)d}\pP\left(\Erw_{0,\omega}\left[\tau_L\right]/L^2\geq
  2\right)\leq \exp\left(-(1/3)(\log L)^2\right).
\end{multline*}
This implies by Borel-Cantelli that
\begin{equation}
\label{times-cormoments-1}
\limsup_{L\rightarrow\infty}\sup_{x: |x| \leq L^k}\sup_{y\in
    V_L(x)}\Erw_{y,\omega}\left[\tau_{V_L(x)}\right]/L^2\leq
  2\quad\pP\hbox{-almost surely}.
\end{equation}
For the rest of the proof, take an environment $\omega$ that 
satisfies~\eqref{times-cormoments-1}.  Assume $m \geq 2$. Then
\begin{align*}
\Erw_{x,\omega}\left[\tau^m_{V_L(x)}\right] &=
\sum_{\ell_1,\dots,\ell_m\geq 0}\Prw_{x,\omega}\left(\tau_{V_L(x)} >
  \ell_1,\dots,\tau_{V_L(x)}>\ell_m\right)\\
&\leq m!\sum_{0\leq \ell_1\leq\dots\leq
  \ell_m}\Prw_{x,\omega}\left(\tau_{V_L(x)} > \ell_m\right).
\end{align*}
By the Markov property, using the case $m=1$ and induction in the last step,
\begin{align*}
\lefteqn{\sum_{0\leq \ell_1\leq\dots\leq
  \ell_m}\Prw_{x,\omega}\left(\tau_{V_L(x)} > \ell_m\right)}\\ 
&= \sum_{0\leq \ell_1\leq\dots\leq
  \ell_{m-1}}\Erw_{x,\omega}\left[\sum_{\ell=0}^\infty\Prw_{X_{\ell_{m-1}},\omega}\left(\tau_{V_L(x)} >
    \ell\right);\, \tau_{V_L(x)}> \ell_{m-1}\right]\\
&\leq \sup_{z\in V_L(x)}\Erw_{z,\omega}\left[\tau_{V_L(x)}\right]\sum_{0\leq \ell_1\leq\dots\leq
  \ell_{m-1}}\Erw_{x,\omega}\left[\tau_{V_L(x)}> \ell_{m-1}\right]\leq 2^mL^{2m},
\end{align*}
if $L = L(\omega)$ is sufficiently large.
\end{proof2}

\section{A quenched invariance principle}
\label{SCLT}
Here we combine the results on the exit distributions from~\cite{BB} and our
results on the mean sojourn times to prove Theorem \ref{thm-clt}, which provides a
functional central limit theorem for the RWRE under the quenched measure.
Let us recall the precise statement.
\begin{quote}
  Assume $\czero(\e)$ for small $\e > 0$, $\cone$ and $\cB$. Then for
  $\pP$-a.e. $\omega\in\Omega$, under $\Prw_{0,\omega}$, the
  $C(\mathbb{R}_+,\mathbb{R}^d)$-valued sequence $X_t^{n}/\sqrt{n}$, $t\geq
  0$, converges in law to a $d$-dimensional Brownian motion with diffusion
  matrix $D^{-1}{\bf \Lambda}$, where $D$ is the constant from
  Theorem~\ref{thm-times}, ${\bf \Lambda}$ is given
  by~\eqref{covariance-matrix-def}, and $X_t^{n}$ is the linear interpolation
  $ X_t^{n}= X_{\lfloor tn\rfloor} + {(tn-\lfloor tn\rfloor)}(X_{\lfloor
    tn\rfloor+1}-X_{\lfloor tn\rfloor}).$
\end{quote}
The statement follows if we show that for each real $T>0$, weak convergence
occurs in $C([0,T],\mathbb{R}^d)$. In order to simplify notation, we
will restrict ourselves to $T=1$, the general case being a simple
generalization of this case.

Let us first give a rough (simplified) idea of our proof. We define the
step size $L_n=(\log n)^{-1}\sqrt{n}$. From Theorem~\ref{thm-times} we
infer that the RWRE should have left about $(\log n)^2/D$ balls of radius
$L_n$ in the first $n$ steps. Proposition~\ref{main-prop-exitmeas} tells us
that for sufficiently large $n$, the exit law from each of those balls is
close to that of a symmetric random walk with nearest neighbor kernel
$p_{L_n}$. For our limit theorem, this will imply that we can replace the
coarse grained RWRE taking steps of size $L_n$, i.e. the RWRE observed at
the successive exit times from balls of radius $L_n$, by the analogous
coarse grained random walk with kernel $p_{L_n}$. For the latter, we apply
the multidimensional Lindeberg-Feller limit theorem. Since we know that the
kernels $p_{L_n}$ converge to $p_\infty$ (see~\eqref{eq-pinfty} and the
comments below Proposition~\ref{main-prop-exitmeas}), we obtain in this way
the stated convergence of the one-dimensional distributions. Since our
estimates on exit measures and exit times are sufficiently uniform in in
the starting point, multidimensional convergence and tightness then follow
from standard arguments.
\begin{figure}[ht]
\begin{center}\parbox{7.5cm}{\includegraphics[width=7cm]{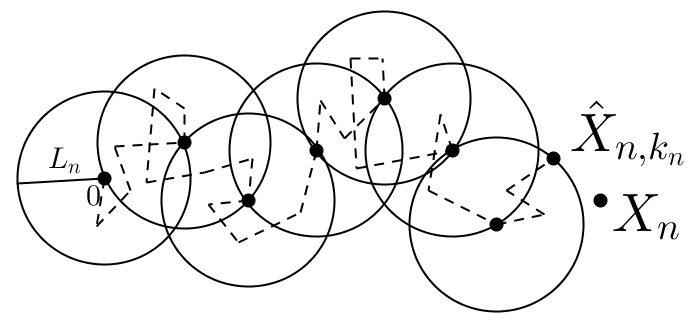}}
\parbox{8cm}{
  \caption{The coarse grained RWRE $\hat{X}_{n,i}$, $i\in\mathbb{N}$, which
    is obtained from observing the RWRE at the successive exit times from
    balls of radius $L_n$. Here $k_n$ denotes the maximal number of such
    balls which are left by the RWRE in the first $n$ steps.}}
\end{center}
\end{figure}

\subsection{Construction of coarse grained random walks on $\mathbb{Z}^d$}
We start with a precise description of the coarse grained random walks. Let
$$L_n=(\log n)^{-1}\sqrt{n}.$$ 
Similarly to the proof of Lemma~\ref{times-keylemma}, given an environment
$\omega\in\Omega$, we introduce a probability space where we can observe
both the random walk with kernel $p_\omega$ and a coarse grained version of
it taking steps of a size between $L_n$ and $2L_n$.  

More specifically, we take a probability space $(\Xi, \mathcal{A}, \pQ)$
that carries a family of i.i.d. random variables $(\xi_{n,i}:
i\in\mathbb{N})$, with $\xi_{n,i}$ distributed according to $\varphi(t)\dt
t$. We then consider the probability space
$((\mathbb{Z}^d)^{\mathbb{N}}\times\Xi,
\pG\otimes\mathcal{A},\Prwt_{x,\omega})$, where $\Prwt_{x,\omega} =
\Prw_{x,\omega}\otimes\pQ$. On this space, $X_k$ denotes again the
projection on the $k$th component of
$\left(\mathbb{Z}^d\right)^{\mathbb{N}}$, so that under $\Prwt_{x,\omega}$,
$X_k$ has the law of a random walk started from $x$ with transition kernel
$p_{\omega}$.

Set $T_{n,0} = 0$, and recursively
for integers $i\in\mathbb{N}$,
\begin{align*}
T_{n,i+1} &= \inf\left\{m> T_{n,i}: X_m\notin V_{\xi_{n,T_{n,i}}\cdot
    L_n}\left(X_{T_{n,i}}\right)\right\}\\
\hat{X}_{n,i} &= X_{T_{n,i}}.
\end{align*}
Under $\Prwt_{x,\omega}$, for fixed $n$, $\hat{X}_{n,i}$ is the coarse
grained Markov chain running with transition
probabilities
$$
Q_{n,\omega}(y,\cdot)=
\frac{1}{L_n}\int_{\mathbb{R}_+}\varphi\left(\frac{t}{L_n}\right)\Pi_{V_t(y),\omega}(y,\cdot)\dt
t
$$
and started from $x$, i.e. $\Prwt_{x,\omega}(\hat{X}_{n,0}=x)=1$. Note that
in contrast to Lemma~\ref{times-keylemma}, the step size of the coarse
grained walk takes values between $L_n$ and $2L_n$ and does not depend on
the current location. We shall suppress the environment $\omega$ in the
notation and write $Q_{n}$ instead of $Q_{n,\omega}$.

We will compare $Q_n$ with the coarse grained (non-random) kernel
$$
q_{n}(y,\cdot) =
\frac{1}{L_n}\int_{\mathbb{R}_+}\varphi\left(\frac{t}{L_n}\right)\pi_{V_t}^{(p_{L_n})}(0,\cdot-y)\dt t,
$$
where the kernel $p_{L_n}$ stems from the assignment~\eqref{kernelpL}.

\subsection{Good events}
\label{clt-prel}
\subsubsection{Good behavior in space}
We shall now introduce an event $A_1$ with $\pP(A_1)=1$ on which the RWRE has a ``good''
behavior in terms of exit distributions. Let
$$D_{L,p,\psi,q}(x) = \left\|\left(\Pi_{V_L(x)}-\pip_{V_L(x)}\right)\piq_{\psi}(x,\cdot)\right\|_1.$$
We require that all smoothed differences of exit measures
$D_{L,p_{L_n},\psi_n,p_{L_n}}(x)$, where $x\in\mathbb{Z}^d$ with $|x|\leq n^3$,
$L\in [L_n,2L_n]$ and $\psi_n\equiv L_n$, are small when $n$ is large.

In this regard, note that Proposition~\ref{main-prop-exitmeas} implies for
large $n$
\begin{align*}
\lefteqn{\pP\left(\sup_{x:|x|\leq
    n^3}\sup_{L_n\leq L\leq 2L_n}D_{L,p_{L_n},\psi_n,p_{L_n}}(x) > (\log L_n)^{-9}\right)}\\
&\leq C n^{4d}\sup_{L_n\leq L\leq 2L_n}\pP\left(D_{L,p_{L_n},\psi_n,p_{L_n}}^\ast > (\log L_n)^{-9}\right)
\leq  \exp\left(-(1/5)(\log n)^2\right).
\end{align*}
An application of Borel-Cantelli then shows that on a set $A_1$ of full $\pP$-measure,
\begin{equation}
\label{eq1}
\limsup_{L\rightarrow\infty}\sup_{x:|x|\leq
    n^3}\sup_{L_n\leq L\leq 2L_n}D_{L,p_{L_n},\psi_n,p_{L_n}}(x) \leq (\log L_n)^{-9}.
\end{equation}

\subsubsection{Good behavior in time}
We next specify an event $A_2$ with $\pP(A_2)=1$ on which we have  
uniform control over mean sojourn times. Let
$$c_\varphi=\int_{1}^2t^2\varphi(t)\dt t.$$
Under our usual conditions, we obtain by Theorem~\ref{thm-times} and dominated convergence, for
$\pP$-almost all $\omega\in\Omega$, $D$ the constant from the theorem,
\begin{equation}
\label{eq-timelimit}
\lim_{n\rightarrow\infty}\left(\inf_{x: |x|\leq n^3}
\Erwt_{x,\omega}\left[T_{n,1}\right]/L_n^2\right) =
\lim_{n\rightarrow\infty}\left(\sup_{x: |x|\leq n^3}
\Erwt_{x,\omega}\left[T_{n,1}\right]/L_n^2\right) =
c_\varphi D.
\end{equation}
Moreover, by Corollary~\ref{cor-times-moments}, for $\pP$-almost all $\omega$,
\begin{equation}
\label{eq-timemoments}
\limsup_{n\rightarrow\infty}\left(\sup_{x: |x|\leq n^3}
\Erwt_{x,\omega}\left[T^2_{n,1}\right]/L_n^4\right) \leq 8.
\end{equation}
We denote by $A_2$ the set of environments of full $\pP$-measure on which
both~\eqref{eq-timelimit} and~\eqref{eq-timemoments} hold true. 

\subsection{A law of large numbers}
Recall Figure $2$. We shall not merely consider $k_n=k_{n,1}$, but more
generally for $t\in[0,1]$
$$
k_{n,t}= k_{n,t}(\omega)=\max\left\{i\in\mathbb{N} : T_{n,i}\leq tn
\right\}.
$$
We shall need a (weak) law of large numbers for $k_{n,t}$ under
$\Prwt_{x,\omega}$, uniformly in $|x|\leq n^2$. In view of~\eqref{eq-timelimit}, it is natural to
expect that $k_{n,t}$ has the same asymptotic behavior as $t\beta_n$, where 
$$
\beta_n= \left\lfloor \frac{n}{c_\varphi DL_n^2}\right\rfloor =
\left\lfloor \frac{(\log n)^2}{c_\varphi D}\right\rfloor 
$$

We first establish a bound on the variance of $T_{n,\ell}$.
\begin{lemma}
\label{clt-lemma-variance}
For $\pP$-almost all environments, 
$$
\sup_{\ell\leq
  2\beta_n}\sup_{|x|\leq n^2}\frac{\textup{Var}_{\Prwt_{x,\omega}}(T_{n,\ell})}{n^2}\rightarrow
0\quad\hbox{as }n\rightarrow\infty,
$$
where $\textup{Var}_{\Prwt_{x,\omega}}$ denotes the variance with respect to $\Prwt_{x,\omega}$.
\end{lemma}
\begin{proof1}
  We can restrict ourselves to $\omega\in A_2$. Define the successive
  sojourn times $\tau_{n,i}= (T_{n,i}-T_{n,i-1})$. Then $T_{n,\ell}=
  \tau_{n,1}+\dots+\tau_{n,\ell}$. Unlike for random walk in a homogeneous
  environment, the variables $\tau_{n,i}$, $i=1,\dots,2\beta_n$, are in
  general not independent under $\Prwt_{0,\omega}$. However, for $i<j$,
  $\tau_{n,j}$ is conditionally independent from $\tau_{n,i}$ given
  $\hat{X}_{n,j-1}$. By the strong Markov property (with the
  same justification as in the proof of Lemma~\ref{times-keylemma}),
  we obtain for $i<j\leq 2\beta_n$ and
  $x\in\mathbb{Z}^d$ with $|x|\leq n^2$,
\begin{align*}
\Erwt_{x,\omega}\left[\tau_{n,i}\tau_{n,j}\right]
&=\Erwt_{x,\omega}\left[\tau_{n,i}\Erwt_{x,\omega}\left[\tau_{n,j}\,\big|\,\hat{X}_{n,j-1},
  \tau_{n,i}\right]\right]\\
&=\Erwt_{x,\omega}\left[\tau_{n,i}\Erwt_{\hat{X}_{n,j-1},\omega}\left[T_{n,1}\right]\right]
\leq \sup_{|y|\leq 2n^2}\Erwt_{y,\omega}\left[T_{n,1}\right]^2.
\end{align*}
In the last step we used that the coarse grained random
can bridge in $2\beta_n$ steps a distance of at most $4\beta_nL_n<n$ and is
therefore well inside $V_{2n^2}$ when started from $V_{n^2}$.
Similarly, we see that
$$\Erwt_{x,\omega}\left[\tau_{n,i}\tau_{n,j}\right] \geq
\inf_{|y|\leq 2n^2}\Erwt_{y,\omega}\left[T_{n,1}\right]^2.$$
For $x$ with $|x|\leq n^2$ and $i,j\leq 2\beta_n$, it also holds that 
$$
\inf_{|y|\leq
2n^2}\Erwt_{y,\omega}\left[T_{n,1}\right]^2 \leq \Erwt_{x,\omega}\left[\tau_{n,i}\right]\Erwt_{x,\omega}\left[\tau_{n,j}\right]
\leq \sup_{|y|\leq
2n^2}\Erwt_{y,\omega}\left[T_{n,1}\right]^2.
$$
Since by definition of the event $A_2$, we have for $\omega\in A_2$
$$
\lim_{n\rightarrow\infty}\left(\inf_{|y|\leq
2n^2}\Erwt_{y,\omega}\left[T_{n,1}\right]^2/L_n^4\right)=
\lim_{n\rightarrow\infty}\left(\sup_{|y|\leq
2n^2}\Erwt_{y,\omega}\left[T_{n,1}\right]^2/L_n^4\right),
$$
we obtain for $i,j\leq 2\beta_n$ and $x$ with $|x|\leq n^2$,
\begin{align*}
  \lefteqn{\left|\Erwt_{x,\omega}\left[\tau_{n,i}\tau_{n,j}\right]-\Erwt_{x,\omega}\left[\tau_{n,i}\right]\Erwt_{x,\omega}\left[\tau_{n,j}\right]\right|}\\
&\leq \sup_{|y|\leq
2n^2}\Erwt_{y,\omega}\left[T_{n,1}\right]^2 - \inf_{|y|\leq
2n^2}\Erwt_{y,\omega}\left[T_{n,1}\right]^2\df \alpha(n)=o(L_n^4)\quad\hbox{for }n\rightarrow\infty.
\end{align*}
Using this for $i\neq j$ and~\eqref{eq-timelimit},~\eqref{eq-timemoments} for the terms with $i=j$,
we conclude that for $n\geq n(\omega)$, $\ell\leq 2\beta_n$,
$$
\sup_{|x|\leq n^2}\textup{Var}_{\Prwt_{x,\omega}}(T_{n,\ell}) \leq  
C\beta_nL_n^4 +C\beta_n^2\alpha(n) = o(n^2).$$
This finishes the proof.
\end{proof1}

We are now in position to prove a weak law of large numbers for $k_{n,t}$.
\begin{lemma}
\label{clt-lemma-k_n}
For $\pP$-almost all environments, for every $t\in[0,1]$ and every $\epsilon>0$, 
$$\sup_{|x|\leq n^2}\Prwt_{x,\omega}\left(\Big|k_{n,t}/\beta_n
    -t\Big|>\epsilon\right) \rightarrow 0\quad\hbox{as }n\rightarrow\infty.$$
\end{lemma}
\begin{proof1}
  We take $\omega\in A_2$ as in the previous lemma. There is nothing to
  show for $t=0$, so assume $t\in(0,1]$.  If the statement would not hold,
  then we could find $\epsilon$, $\epsilon'>0$ such that
\begin{align}
\label{clt-k_n1}
\sup_{|x|\leq n^2}\Prwt_{x,\omega}\left(k_{n,t}<
  (t-\epsilon)\beta_n\right)&>\epsilon'\quad\hbox{ infinitely often, or}\\
\label{clt-k_n2}
\sup_{|x|\leq n^2}\Prwt_{x,\omega}\left(k_{n,t}>
  (t+\epsilon)\beta_n\right)&>\epsilon'\quad\hbox{ infinitely often}.
\end{align}
Let us first assume \eqref{clt-k_n1}. Then, with $i_n= \lceil(t-\epsilon)\beta_n\rceil$, by definition
$$
\sup_{|x|\leq n^2}\Prwt_{x,\omega}\left(T_{n,i_n}>tn\right)>\epsilon'\quad\hbox{ infinitely often}.
$$
Next note that by~\eqref{eq-timelimit}, by linearity of the expectation and
the fact that $2i_nL_n<n$, 
$$
0\leq \frac{\sup_{|x|\leq n^2}\Erwt_{x,\omega}\left[T_{n,i_n}\right]}{tn}\leq
\frac{i_n}{tn}\sup_{y:|y|\leq 2n^2}\Erwt_{y,\omega}\left[T_{n,1}\right]\leq
1-\epsilon/2\quad\hbox{for } n\geq n_0(\omega).$$ 
Chebycheff's inequality then shows that if $n\geq n_0(\omega)$ and $x$ with
$|x|\leq n^2$,
\begin{align*}
\Prwt_{x,\omega}\left(T_{n,i_n}>tn\right)&=
\Prwt_{x,\omega}\left(T_{n,i_n}-\Erwt_{x,\omega}\left[T_{n,i_n}\right]> tn-\Erwt_{x,\omega}\left[T_{n,i_n}\right]\right)\\
&\leq
\frac{1}{(tn)^2}\left(1-\Erwt_{x,\omega}\left[T_{n,i_n}\right]/(tn)\right)^{-2}\,\hbox{Var}_{\Prwt_{x,\omega}}\left(T_{n,i_n}\right)\\
&\leq \frac{4}{(\epsilon tn)^2}\sup_{|y|\leq
    n^2}\hbox{Var}_{\Prwt_{y,\omega}}\left(T_{n,i_n}\right).
\end{align*}
The right hand side converges to zero by Lemma~\ref{clt-lemma-variance}. 
This contradicts~\eqref{clt-k_n1}. 

Now assume~\eqref{clt-k_n2}. We argue similarly. First, with $i_n=
\lfloor(t+\epsilon)\beta_n\rfloor$ by definition 
$$
\sup_{|x|\leq n^2}\Prwt_{x,\omega}\left(T_{n,i_n}<tn\right)>\epsilon'\quad\hbox{ infinitely often}.
$$
Since for $\omega\in A_2$ and large $n\geq n_0(\omega)$,
$$\frac{\inf_{|x|\leq n^2}\Erwt_{x,\omega}\left[T_{n,i_n}\right]}{tn}\geq \frac{i_n}{tn}
\inf_{y:|y|\leq 2n^2}\Erwt_{y,\omega}\left[T_{n,1}\right]\geq
1+\epsilon/2,$$
we obtain for large $n\geq n_0(\omega)$ and $x$ with $|x|\leq n^2$,
\begin{align*}
\Prwt_{x,\omega}\left(T_{n,i_n}<tn\right) &=
\Prwt_{x,\omega}\left(\Erwt_{x,\omega}\left[T_{n,i_n}\right]-T_{n,i_n}>
  \Erwt_{x,\omega}\left[T_{n,i_n}\right]-tn\right)\\
& \leq \frac{4}{(\epsilon tn)^2}\sup_{|y|\leq n^2}\hbox{Var}_{\Prwt_{y,\omega}}\left(T_{n,i_n}\right)/n^2\rightarrow
0\quad\hbox{as }n\rightarrow\infty.
\end{align*}
Therefore, neither~\eqref{clt-k_n1} nor~\eqref{clt-k_n2} can hold, and the
proof of the lemma is complete.
\end{proof1}

\subsection{Proof of Theorem~\ref{thm-clt}}
We turn to the proof of Theorem~\ref{thm-clt}. Recall our notation
introduced above. Since the subscript
$n$ already appears in both $k_{n,t}$ and $\beta_n$, we may safely write
$$\hat{X}_{k_{n,t}}\hbox{ instead of }\hat{X}_{n,k_{n,t}},\quad
\hat{X}_{\lfloor t\beta_n\rfloor}\hbox{ instead of }\hat{X}_{n, \lfloor
  t\beta_n\rfloor}.$$ 

Since both $A_1$ and $A_2$ have full $\pP$-measure, we can restrict
ourselves to $\omega\in A_1\cap A_2$.  We first prove one-dimensional
convergence, uniformly in the starting point $x$ with $|x|\leq n^2$. This
will easily imply multidimensional convergence and tightness.

\subsubsection{One-dimensional convergence}
\begin{proposition}
\label{prop1d}
For $\pP$-almost all environments, for each $t\in[0,1]$ and
$u\in\mathbb{R}$, $$ \sup_{|x|\leq
  n^2}\left|\Prw_{x,\omega}\left(\left(X_t^{n}-x\right)/\sqrt{n}>u\right)-\nP\left(\mathcal{N}(0,tD^{-1}{\bf \Lambda})>u\right)\right|\rightarrow
0\quad\hbox{as }n\rightarrow\infty,
$$
where $\mathcal{N}(0,A)$ denotes a $d$-dimensional centered Gaussian with
covariance matrix $A$.
\end{proposition}
\begin{proof1}
Let $t\in[0,1]$. We write
$$X_t^n= \hat{X}_{\lfloor t\beta_n\rfloor}+(X_t^n-\hat{X}_{k_{n,t}}) + (\hat{X}_{k_{n,t}}-
\hat{X}_{\lfloor t\beta_n\rfloor}).$$ Since by definition of the random sequence $k_{n,t}$, one has
$$
\big|X_t^n -\hat{X}_{k_{n,t}}\big|\leq 1 + \big|X_{\lfloor tn\rfloor}
-\hat{X}_{k_{n,t}}\big|\leq 3L_n=o\left(\sqrt{n}\right),
$$ 
our claim follows from the
following two convergences when $n\rightarrow\infty$.
\begin{enumerate}
\item For each $u\in\mathbb{R}$,
    $$\sup_{|x|\leq
      2n^2}\left|\Prwt_{x,\omega}\left(\left(\hat{X}_{\lfloor
            t\beta_n\rfloor}-x\right)/\sqrt{n}>u\right)-\nP\left(\mathcal{N}(0,tD^{-1}{\bf \Lambda})>u\right)\right|\rightarrow
    0.$$
\item For each $\epsilon >0$, $\sup_{|x|\leq
    n^2}\Prwt_{x,\omega}\left(\big|\hat{X}_{k_{n,t}}-\hat{X}_{\lfloor
      t\beta_n\rfloor}\big|/\sqrt{n}> \epsilon\right)\rightarrow 0$.
\end{enumerate}
We first prove (i). For notational simplicity, we restrict ourselves to the
case $t=1$; the general case $t\in[0,1]$ follows exactly the same lines,
with $\beta_n$ replaced everywhere by $\lfloor t\beta_n\rfloor$. For later
use, it will be helpful to consider here the supremum over $x$ bounded by
$2n^2$ instead of $n^2$. We let $(Z_{n,i})_{i=0,\dots,n}$ be an
i.i.d. sequence of random vectors distributed according to $q_n(0,\cdot)$,
independently of the RWRE. Since $|Z_{n,i}|\leq 2L_n =o(\sqrt{n})$, it
suffices to show the statement for $\hat{X}_{\beta_n}$ inside the
probability replaced by $\hat{X}_{\beta_n} + Z_{n,0}$ (tacitly assuming
that $\hat{X}_{\beta_n}$ under $\Prwt_{x,\omega}$ and $Z_{n,0}$ are defined
on the same probability space, whose probability measure we again denote by
$\Prwt_{x,\omega}$).  Now let $\hat{Y}_{i}= Z_{n,1}+\dots+Z_{n,i}$. Since
$\hat{X}_{\beta_n}+Z_{n,0}$ has law $(Q_n)^{\beta_n}q_n(x,\cdot)$ under
$\Prwt_{x,\omega}$, and $x + \hat{Y}_{i}$ has law $(q_n)^{i}(x,\cdot)$, we
get
\begin{multline*}
\sup_{|x|\leq 2n^2}\left|\Prwt_{x,\omega}\left(\left(\hat{X}_{\beta_n}+Z_{n,0}-x\right)/\sqrt{n}>u\right)-
\nP\left((x + \hat{Y}_{\beta_n +1}-x)/\sqrt{n}>u\right)\right|\\
 \leq \sup_{|x|\leq 2n^2}\left\|\left((Q_n)^{\beta_n}-(q_n)^{\beta_n}\right)q_n(x,\cdot)\right\|_1.
\end{multline*}
For $\omega\in A_1$, we obtain by iteration, uniformly in $x$ with $|x|\leq 2n^2$,
\begin{align*}
\lefteqn{\left\|\left((Q_n)^{\beta_n}-(q_n)^{\beta_n}\right)q_n(x,\cdot)\right\|_1}\\
&\leq\left\|(Q_n)^{\beta_n-1}\left(Q_n-q_n\right)q_n(x,\cdot)\right\|_1
+\left\|\left((Q_n)^{\beta_n-1}-(q_n)^{\beta_n-1}\right)q_n^2(x,\cdot)\right\|_1\\
&\leq  \sup_{|x|\leq 3n^2}\left\|\left(Q_n-q_n\right)q_n(x,\cdot)\right\|_1 + \sup_{|x|\leq
  2n^2}\left\|\left((Q_n)^{\beta_n-1}-(q_n)^{\beta_n-1}\right)q_n(x,\cdot)\right\|_1\\
&\leq \beta_n(\log L_n)^{-9}\rightarrow 0\quad\hbox{as }n\rightarrow\infty.
\end{align*}
It remains to show that $\hat{Y}_{\beta_n}/\sqrt{n}$ converges in
distribution to a $d$-dimensional centered Gaussian vector with covariance
matrix $D^{-1}{\bf \Lambda}$. This will be a consequence of the following
multidimensional version of the Lindeberg-Feller theorem.
\begin{proposition}
\label{lindfel-prop}
  Let $W_{m,\ell}$, $1\leq \ell\leq m$, be centered and independent 
  $\mathbb{R}^d$-valued random vectors. Put ${\bf\Sigma}_{m,\ell} = (\sigma_{m,\ell}^{(ij)})_{i,j=1,\dots, d}$, 
  where $\sigma_{m,\ell}^{(ij)}=
  \nE\left[W_{m,\ell}^{(i)}W_{m,\ell}^{(j)}\right]$ and
  $W_{m,\ell}^{(i)}$ is the $i$th component of $W_{m,\ell}$. If
  for $m\rightarrow \infty$,
\begin{enumerate}
\item[(a)] $\sum_{\ell=1}^m{\bf\Sigma}_{m,\ell}\rightarrow {\bf\Sigma}$, 
\item[(b)] for each ${\bf v}\in\mathbb{R}^d$ and each $\epsilon> 0$, 
$\sum_{\ell=1}^m\nE\left[\left|{\bf v}\cdot W_{m,\ell}\right|^2;\, \left|{\bf v}\cdot W_{m,\ell}\right|
  > \epsilon\right]\rightarrow 0$,
\end{enumerate}
then $W_{m,1}+\dots + W_{m,m}$ converges in distribution as $m\rightarrow\infty$
to a $d$-dimensional Gaussian random vector with mean zero and
covariance matrix ${\bf\Sigma}$.
\end{proposition}

\begin{proof2}{\bf of Proposition~\ref{lindfel-prop}:}
  By the Cram\'er-Wold device, it suffices to show that for fixed
  ${\bf v}\in\mathbb{R}^d$, ${\bf v}\cdot (W_{m,1}+\dots+W_{m,m})$ converges
  in distribution to a Gaussian random variable with mean zero and variance
  ${\bf v}^T{\bf\Sigma}{\bf v}$. Under (a) and (b), this follows immediately
  from the classical one-dimensional Lindeberg-Feller theorem.
\end{proof2}
We now finish the proof of (i). Recall that
$\hat{Y}_{\beta_n}= Z_{n,1}+\dots+Z_{n,\beta_n}$, where the $Z_{n,\ell}$
are independent random vectors with law $q_n(0,\cdot)$. Since the
underlying one-step transition kernel $p_{L_n}$ is symmetric, the
$Z_{n,\ell}$ are centered. Moreover, denoting by
$Z_{n,\ell}^{(i)}$ the $i$th component of $Z_{n,\ell}$,
$$
\nE\left[Z_{n,\ell}^{(i)}Z_{n,\ell}^{(j)}\right]=0\quad\hbox{for
}i\neq j,\quad i,j=1,\dots,d.
$$
For $i=j$, we obtain by definition
\begin{align}
\label{eiequalsj}
\lefteqn{\frac{1}{n}\left(\nE\left[\left(Z_{n,1}^{(i)}\right)^2\right]+\dots+\nE\left[\left(Z_{n,\beta_n}^{(i)}\right)^2\right]\right)}\nonumber\\
&= \frac{\beta_n}{n}\sum_{y\in\mathbb{Z}^d}q_n(0,y)y_i^{2}
=\frac{\beta_n}{n}\int_1^2\varphi(s)\sum_{y\in\mathbb{Z}^d}\pi_{V_{sL_n}}^{(p_{L_n})}(0,y)y_i^{2}\dt s\nonumber\\
&=\frac{\beta_n}{n}\int_1^2(L_ns)^2\varphi(s)\sum_{y\in\mathbb{Z}^d}\pi_{V_{sL_n}}^{(p_{L_n})}(0,y)(y_i/sL_n)^{2}\dt s.
\end{align}
We next recall that~\cite[Lemma 2.1]{BB} shows how to recover the kernel
$p_{L_n}$ out of the exit measure $\pi_{V_{sL_n}}^{(p_{L_n})}$, namely
$$
2p_{L_n}=\sum_{y\in\mathbb{Z}^d}\pi_{V_{sL_n}}^{(p_{L_n})}(0,y)(y_i/sL_n)^{2}
+ O(L_n^{-1}).
$$
Replacing $\beta_n$ by its value, we therefore deduce from~\eqref{eiequalsj}
that
\begin{align*}
\lefteqn{\frac{1}{n}\left(\nE\left[\left(Z_{n,1}^{(i)}\right)^2\right]+\dots+\nE\left[\left(Z_{n,\beta_n}^{(i)}\right)^2\right]\right)}\\
&=\frac{L_n^2\beta_nc_\varphi}{n}2p_{L_n}(e_i) +
O\left(L_n^{-1}\right)
=\frac{2p_{L_n}(e_i)}{D} + O\left((\log n)^{-2}\right).
\end{align*}
Since $p_{L_n}(e_i)\rightarrow p_\infty(e_i)$ as $n\rightarrow\infty$, we
obtain with $ W_{\beta_n,\ell} = Z_{n,\ell}/\sqrt{n}$,
$\ell=1,\dots,\beta_n$, in the notation of Proposition~\ref{lindfel-prop},
$$\sum_{\ell=1}^{\beta_n}{\bf\Sigma}_{\beta_n,\ell}\rightarrow
D^{-1}\left(2p_\infty(e_i)\delta_i(j)\right)_{i,j=1}^d=D^{-1}{\bf \Lambda}\quad\hbox{as
}n\rightarrow\infty.$$ Since $W_{\beta_n,\ell}\leq 2L_n/\sqrt{n}\leq 2(\log
n)^{-1}$, point (b) of Proposition~\ref{lindfel-prop} is trivially
fulfilled. Applying this proposition finally shows that $
\hat{Y}_{\beta_n}/\sqrt{n}=W_{ \beta_n,1}+\dots+W_{\beta_n,\beta_n}$
converges in distribution to a $d$-dimensional centered Gaussian random
vector with covariance matrix $D^{-1}{\bf \Lambda}$. This finishes the proof of
(i).

It remains to prove (ii). In view of Lemma \ref{clt-lemma-k_n}, it suffices
to show that for each $\epsilon>0$,
$$
\lim_{\theta\downarrow 0}\lim_{n\rightarrow\infty}\sup_{|x|\leq n^2}\Prwt_{x,\omega}\left(\big|\hat{X}_{k_{n,t}}-\hat{X}_{\lfloor
      t\beta_n\rfloor}\big|>\epsilon\sqrt{n};\,\big|k_{n,t}-\lfloor
    t\beta_n\rfloor \big|\, < \theta\beta_n\right) = 0.
$$ 
Fix $\epsilon>0$, $\theta>0$. Define the set of integers
$$A_{n}=\{\lfloor
t\beta_n\rfloor - \lceil\theta\beta_n\rceil,\dots,\lfloor t\beta_n\rfloor
+ \lceil\theta\beta_n\rceil\},$$ and let $\underline{\ell_n}=\lfloor
t\beta_n\rfloor - \lceil\theta\beta_n\rceil$.  Then
\begin{align*}
\lefteqn{\Prwt_{x,\omega}\left(\big|\hat{X}_{k_{n,t}}-\hat{X}_{\lfloor
      t\beta_n\rfloor}\big|>\epsilon\sqrt{n};\,\big|k_{n,t}-\lfloor
    t\beta_n\rfloor \big|\, < \theta\beta_n\right)}\\
  &\leq \Prwt_{x,\omega}\left(\max_{\ell\in A_n}\big|\hat{X}_{n,\ell}-\hat{X}_{\lfloor
      t\beta_n\rfloor}\big|>\epsilon\sqrt{n}\right)\\
&\leq \Prwt_{x,\omega}\left(\max_{\ell\in
    A_n}\big|\hat{X}_{n,\ell}-\hat{X}_{n,\underline{\ell_n}}\big|>(\epsilon/2)\sqrt{n}\right)
  + \Prwt_{x,\omega}\left(\big|\hat{X}_{n,\underline{\ell_n}}-\hat{X}_{\lfloor
      t\beta_n\rfloor}\big|>(\epsilon/2)\sqrt{n}\right). 
\end{align*}
We only consider the first probability in the last display; the second one
is treated in a similar (but simpler) way. We first remark that after
$\underline{\ell_n}$ steps, the coarse grained RWRE with transition kernel
$Q_n$ starting in $V_{n^2}$ is still within $V_{2n^2}$. Therefore, by the
Markov property, for $x$ with $|x|\leq n^2$,
\begin{align}
\label{eq-1d-1}
\lefteqn{\Prwt_{x,\omega}\left(\max_{\ell\in A_n}\big|\hat{X}_{n,\ell}-\hat{X}_{n,\underline{\ell_n}}\big|>(\epsilon/2)\sqrt{n}\right)}\nonumber\\
&\leq \sup_{|y|\leq 2n^2} \Prwt_{y,\omega}\left(\max_{\ell\leq
    2\lceil\theta\beta_n\rceil}\big|\hat{X}_{n,\ell}-y\big|>(\epsilon/2)\sqrt{n}\right).
\end{align}
For estimating~\eqref{eq-1d-1}, we follow a strategy similar to
Billingsley~\cite[Theorem 9.1]{BIL1}. Put
$$E_\ell=\left\{\max_{j<\ell}\big|\hat{X}_{n,j}-\hat{X}_{n,0}\big|<(\epsilon/2)\sqrt{n}\leq
  \big|\hat{X}_{n,\ell}-\hat{X}_{n,0}\big|\right\}.$$ Then
\begin{multline*}
\Prwt_{y,\omega}\left(\max_{\ell\leq  2\lceil\theta\beta_n\rceil}\big|\hat{X}_{n,\ell}-y\big|>
    (\epsilon/2)\sqrt{n}\right)\leq \Prwt_{y,\omega}\left(\big|\hat{X}_{n,2\lceil\theta\beta_n\rceil}-y\big|\geq
    (\epsilon/4)\sqrt{n}\right)\\ 
+ \sum_{\ell=1}^{2\lceil\theta\beta_n\rceil-1}\Prwt_{y,\omega}\left(\big|\hat{X}_{n,
      2\lceil\theta\beta_n\rceil}-\hat{X}_{n,\ell}\big| \geq (\epsilon/4)\sqrt{n};\,E_\ell\right).
\end{multline*}
Concerning the first probability on the right, we already know from  
(i) that for $\theta < 1/2$,
$$ \sup_{|y|\leq
      2n^2}\Prwt_{y,\omega}\left(\big|\hat{X}_{n,2\lceil\theta\beta_n\rceil}-y\big|\geq
    (\epsilon/4)\sqrt{n}\right)\rightarrow\nP\left(\big|\mathcal{N}(0,2\theta
    D^{-1}{\bf \Lambda})\big|\geq\epsilon/4\right)\quad\hbox{as
    }n\rightarrow\infty.$$
For fixed $\epsilon$, the right side converges to zero as
$\theta\downarrow 0$ by Chebycheff's inequality. 
For the sum over the probabilities in the above display, we stress that
the increments of the coarse grained walk $\hat{X}_{n,\ell}$ are neither
independent nor stationary under $\Prwt_{y,\omega}$. But we have by the Markov property at time $\ell$, for
$|y|\leq 2n^2$,
\begin{align}
\label{eq-1d-2}
\lefteqn{\sum_{\ell=1}^{2\lceil\theta\beta_n\rceil-1}\Prwt_{y,\omega}\left(\big|\hat{X}_{n,
        2\lceil\theta\beta_n\rceil}-\hat{X}_{n,\ell}\big| \geq (\epsilon/4)\sqrt{n};\,E_\ell\right)}\nonumber\\
&\leq \sum_{\ell=1}^{2\lceil\theta\beta_n\rceil-1}
\Prwt_{y,\omega}(E_\ell)\sup_{|z|\leq
  3n^2}\Prwt_{z,\omega}\left(\big|\hat{X}_{n,2\lceil\theta\beta_n\rceil-\ell}-z\big| \geq
  (\epsilon/4)\sqrt{n}\right).
\end{align}
Similar to the proof of (i), we estimate for $\ell=1,\dots,2\lceil\theta\beta_n\rceil-1$
\begin{multline*}
\sup_{|z|\leq 3n^2}\Prwt_{z,\omega}\left(\big|\hat{X}_{n,\ell}-z \big| \geq
    (\epsilon/4)\sqrt{n}\right) 
\leq  \sup_{|z|\leq 3n^2}(Q_n)^{\ell}q_n\left(z,\mathbb{Z}^d\backslash V_{(\epsilon/8)\sqrt{n}}(z)\right)\\
\leq \sup_{|z|\leq 3n^2}\left\|\left((Q_n)^{\ell}-(q_n)^{\ell}\right)q_n(z,\cdot)\right\|_1 + 
(q_n)^{\ell+1}\left(0,\mathbb{Z}^d\backslash
  V_{(\epsilon/8)\sqrt{n}}\right).
\end{multline*}
For environments $\omega\in A_1$, the first summand is estimated by $\ell
(\log L_n)^{-9}$ as in the proof of (i). For the expression involving
$q_n$, we use the following standard large deviation estimate (a proof is 
for example given in~\cite[Proof of Lemma 7.5]{BB}): There
exist constants $C_1$, $c_1$ depending only on the dimension such that
$$
(q_n)^{\ell}\left(0,\mathbb{Z}^d\backslash V_{r}\right)\leq
C_1\exp\left(-c_1r^2/(\ell L_n^2)\right),\quad r>0,\,\ell\in\mathbb{N}.
$$
In our setting, we obtain 
$$
(q_n)^{\ell}\left(0,\mathbb{Z}^d\backslash
  V_{(\epsilon/8)\sqrt{n}}\right)\leq C
\exp\left(-c\epsilon^2/\theta\right)\quad\hbox{uniformly in }1\leq\ell\leq
2\lceil\theta\beta_n\rceil.
$$
Back to~\eqref{eq-1d-2}, the fact that the $E_i$'s are disjoint leads to
\begin{align*}
  \lefteqn{\sum_{\ell=1}^{2\lceil\theta\beta_n\rceil-1}\Prwt_{y,\omega}\left(\big|\hat{X}_{n,
        2\lceil\theta\beta_n\rceil}-\hat{X}_{n,\ell}\big| \geq (\epsilon/4)\sqrt{n};\,E_\ell\right)}\\
  &\leq \beta_n^2(\log L_n)^{-9} +
  \sum_{\ell=1}^{2\lceil\theta\beta_n\rceil-1}\Prwt_{y,\omega}\left(E_\ell\right)(q_n)^{2\lceil\theta\beta_n\rceil+1-\ell}\left(0,\mathbb{Z}^d\backslash
    V_{(\epsilon/8)\sqrt{n}}\right)\\
  &\leq o(1) + C\exp\left(-c\epsilon^2/\theta\right),
\end{align*}
everything uniformly in $|y|\leq 2n^2$. The last expression converges to zero as
$\theta\downarrow 0$. This concludes the proof of (ii) and hence of the
one-dimensional convergence.
\end{proof1}

\subsubsection{Convergence of finite-dimensional distributions}
In order to prove convergence of the two-dimensional distributions under
$\Prw_{0,\omega}$, we have to show that for $0\leq t_1<t_2\leq
1$ and $u_1,u_2\in\mathbb{R}$, as $n\rightarrow\infty$,
\begin{multline}
\Big|\Prw_{0,\omega}\left(X_{t_1}^n/\sqrt{n} > u_1,\,(X_{
    t_2}^n-X_{t_1}^n)/\sqrt{n}>u_2\right)\\
 - \nP\left(\mathcal{N}(0,t_1D^{-1}{\bf \Lambda})>u_1\right)\nP\left(\mathcal{N}(0,(t_2-t_1)D^{-1}{\bf \Lambda})>u_2\right)\Big|\rightarrow 0.
\end{multline} 
This follows easily from our uniform one-dimensional convergence. First, we
may replace $X_{t_1}^n$ by $X_{\lfloor t_1n\rfloor}$ and $X_{t_2}^n$ by
$X_{\lfloor t_2n\rfloor}$, since their difference is bounded by one. Then,
by the Markov property
\begin{align*}
\lefteqn{\Prw_{0,\omega}\left(X_{\lfloor t_1n\rfloor}/\sqrt{n} > u_1,\,(X_{\lfloor
    t_2n\rfloor}-X_{\lfloor t_1n\rfloor})/\sqrt{n}>u_2\right)}\\
&=\Prw_{0,\omega}\left(X_{\lfloor t_1n\rfloor}/\sqrt{n} > u_1,\,
  \Prw_{X_{\lfloor t_1n\rfloor},\omega}\left(\left(X_{\lfloor
    t_2n\rfloor-\lfloor
    t_1n\rfloor}-X_0\right)/\sqrt{n}>u_2\right)\right)\\
&\leq \Prw_{0,\omega}\left(X_{\lfloor t_1n\rfloor}/\sqrt{n} >
  u_1\right)\sup_{|x|\leq n}\Prw_{x,\omega}\left(\left(X_{\lfloor
    t_2n\rfloor-\lfloor
    t_1n\rfloor}-x\right)/\sqrt{n}>u_2\right).
\end{align*}
The product of the two probabilities converges by Proposition~\ref{prop1d}
towards
\begin{equation}
\label{ePP}\nP\left(\mathcal{N}(0,t_1D^{-1}{\bf \Lambda})>u_1\right)\nP\left(\mathcal{N}(0,(t_2-t_1)D^{-1}{\bf \Lambda})>u_2\right).
\end{equation}
For the lower bound,
\begin{align*}
\lefteqn{\Prw_{0,\omega}\left(X_{\lfloor t_1n\rfloor}/\sqrt{n} > u_1,\,(X_{\lfloor
    t_2n\rfloor}-X_{\lfloor t_1n\rfloor})/\sqrt{n}>u_2\right)}\\
&\geq \Prw_{0,\omega}\left(X_{\lfloor t_1n\rfloor}/\sqrt{n} >
  u_1\right)\inf_{|x|\leq n}\Prw_{x,\omega}\left(\left(X_{\lfloor
    t_2n\rfloor-\lfloor
    t_1n\rfloor}-x\right)/\sqrt{n}>u_2\right),
\end{align*}
and the right hand side converges again towards the product
in~\eqref{ePP}. This proves convergence of the two-dimensional
distributions under $\Prw_{0,\omega}$. The general case of
finite-dimensional convergence is obtained similarly.

\subsubsection{Tightness}
The sequence of $\Prw_{0,\omega}$-laws of
$(X_t^n/\sqrt{n} : 0\leq t\leq 1)$ is tight, if the following Condition
{\bf T} holds true.
\begin{quote}
  For each $\epsilon >0$ there exist a $\lambda>1$ and an integer $n_0$
  such that, if $n\geq n_0$,
$$
\Prw_{0,\omega}\left(\max_{\ell\leq n}\big|X_{k+\ell}-X_k\big| \geq
  \lambda\sqrt{n}\right) \leq \frac{\epsilon}{\lambda^2}\quad\hbox{for all }k\leq n\lambda^2/\epsilon. 
$$
\end{quote}
See~\cite[Theorem 8.4]{BIL1} for a proof of this standard criterion. 

Let us now show that Condition {\bf T} is indeed satisfied in our
setting. First, by the Markov property at time $k$,
$$
\Prw_{0,\omega}\left(\max_{\ell\leq n}\big|X_{k+\ell}-X_k\big| \geq
  \lambda\sqrt{n}\right)\leq \sup_{|x|\leq k}
\Prw_{x,\omega}\left(\max_{\ell\leq n}\big|X_{\ell}-x\big| \geq \lambda\sqrt{n}\right).
$$
The random walk $X_k$ under $\Prw_{x,\omega}$ has the same law as the first
coordinate process on $(\mathbb{Z}^d)^{\mathbb{N}}\times\Xi$ under
$\Prwt_{x,\omega}$, which we also denote by $X_k$ (see the beginning of
Section~\ref{SCLT}).  We shall now consider the latter under
$\Prwt_{x,\omega}$. We recall that $k_{n,1}=k_{n,1}(\omega)$ counts the
number of steps the coarse grained walk performs up to time $n$. Now we
have
\begin{align*}
\lefteqn{\Prw_{x,\omega}\left(\max_{\ell\leq n}\big|X_{\ell}-x\big|\geq \lambda\sqrt{n}\right)}\\
&\leq\Prwt_{x,\omega}\left(\max_{\ell\leq n}\big|X_{\ell}-x\big|\geq \lambda\sqrt{n};\,
  k_{n,1}\leq 2\beta_n\right) + \Prwt_{x,\omega}\left(k_{n,1}>
  2\beta_n\right).
\end{align*}
The second probability on the right converges to zero as $n$ tends to
infinity by Lemma~\ref{clt-lemma-k_n}, uniformly in starting points $x$
with $|x|\leq n^2$.  For the first probability, we find on the
event $\{k_{n,1}\leq 2\beta_n\}$ for each $j\leq n$ an $\ell\leq 2\beta_n$
such that $|X_j-\hat{X}_{n,\ell}| \leq 2L_n$. We therefore obtain for large
$n$
$$
\Prwt_{x,\omega}\left(\max_{\ell\leq n}\big|X_{\ell}-x\big|\geq \lambda\sqrt{n};\,
  k_{n,1}\leq 2\beta_n\right)\leq \Prwt_{x,\omega}\left(\max_{\ell\leq 2\beta_n}\big|\hat{X}_{n,\ell}-x\big|\geq (\lambda/2)\sqrt{n}\right).
$$
For bounding this last probability, we can follow the same steps as
for estimating~\eqref{eq-1d-1}. Leaving out the details, we arrive at
$$
\sup_{|x|\leq n^2}\Prwt_{x,\omega}\left(\max_{\ell\leq
    2\beta_n}\big|\hat{X}_{n,\ell}-x\big|\geq
  (\lambda/2)\sqrt{n}\right)\leq \frac{C}{\lambda^3} + C\exp\left(-c\lambda^2\right)\leq  \frac{\epsilon}{\lambda^2},
$$
provided $\lambda=\lambda(d,\epsilon)$ is large enough. This proves that Condition {\bf T} is
satisfied. Therefore, the sequence of $\Prw_{0,\omega}$-laws of
$(X_t^n/\sqrt{n} : 0\leq t\leq 1)$ is tight, which concludes also the proof of Theorem~\ref{thm-clt}.
\subsection*{Acknowledgments}
I am indebted to Erwin Bolthausen for constant support and
advice. Furthermore I would like to thank Jean-Christophe Mourrat and Ofer
Zeitouni for helpful discussions.

\end{document}